\documentclass[12pt]{article}

\usepackage{indentfirst}
\usepackage{amsfonts}
\usepackage{amsmath}
\usepackage{amssymb}
\usepackage{amsbsy, amsthm}
\usepackage[margin=0.9in]{geometry}

\newtheorem{dfn}{Definition}[section]
\newtheorem{theorem}[dfn]{Theorem}
\newtheorem{lemma}[dfn]{Lemma}
\newtheorem{proposition}[dfn]{Proposition}

\newtheorem{corollary}[dfn]{Corollary}
\newenvironment{pf}{\noindent{\bf Proof.}}
{\enspace\vrule height5pt depth0pt width5pt}

\def\X {{\mathcal X}}

\def\E {{\mathcal E}}

\def\Se {{\mathcal S}}
\def\F {{\mathcal F}}

\def\D {{\mathcal D}}
\def\Q {{\mathcal Q}}

\def\I {{\mathcal I}}

\def\adh {{\rm adh}}
\def\C {{\mathcal C}}

\begin{document}
\title{A global decomposition theorem for excluding immersions in graphs with no edge-cut of order three}
\author{Chun-Hung Liu\thanks{chliu@math.tamu.edu. Partially supported by NSF under Grant No.~DMS-1929851 and DMS-1954054.} \\
\small Department of Mathematics, \\
\small Texas A\&M University,\\
\small College Station, TX 77843-3368, USA}

\maketitle

\begin{abstract}
\noindent
A graph $G$ contains another graph $H$ as an immersion if $H$ can be obtained from a subgraph of $G$ by splitting off edges and removing isolated vertices.
There is an obvious necessary degree condition for the immersion containment: if $G$ contains $H$ as an immersion, then for every integer $k$, the number of vertices of degree at least $k$ in $G$ is at least the number of vertices of degree at least $k$ in $H$.
In this paper, we prove that this obvious necessary condition is ``nearly'' sufficient for graphs with no edge-cut of order 3: for every graph $H$, every $H$-immersion free graph with no edge-cut of order 3 can be obtained by an edge-sum of graphs, where each of the summands is obtained from a graph violating the obvious degree condition by adding a bounded number of edges. 
The condition for having no edge-cut of order 3 is necessary.
A simple application of this theorem shows that for every graph $H$ of maximum degree $d \geq 4$, there exists an integer $c$ such that for every positive integer $m$, there are at most $c^m$ unlabelled $d$-edge-connected $H$-immersion free $m$-edge graphs with no isolated vertex, while there are superexponentially many unlabelled $(d-1)$-edge-connected $H$-immersion free $m$-edge graphs with no isolated vertex.
Our structure theorem will be applied in a forthcoming paper about determining the clustered chromatic number of the class of $H$-immersion free graphs.
\end{abstract}

\section{Introduction}
All {\it graphs} in this paper are finite and allowed to have loops and parallel edges.
For graphs $G$ and $H$, we say that $G$ contains an {\it $H$-immersion} (or $G$ contains $H$ as an {\it immersion}) if there exist functions $\pi_V$ and $\pi_E$ such that 
	\begin{itemize}
		\item $\pi_V$ is an injection from $V(H)$ to $V(G)$,
		\item $\pi_E$ maps each edge of $H$ to a subgraph of $G$ such that for each $e \in E(H)$, if $e$ has distinct ends $x,y$, then $\pi_E(e)$ is a path in $G$ with ends $\pi_V(x)$ and $\pi_V(y)$, and if $e$ is a loop with end $v$, then $\pi_E(e)$ is a cycle containing $\pi_V(v)$, and
		\item if $e_1,e_2$ are distinct edges of $H$, then $\pi_E(e_1)$ and $\pi_E(e_2)$ are edge-disjoint.
	\end{itemize}
The immersion containment is closely related to the minor and topological minor containments.
A graph $G$ contains another graph $H$ as a {\it topological minor} if some subgraph of $G$ is isomorphic to a subdivision of $H$.
And $G$ contains $H$ as a {\it minor} if $H$ is isomorphic to a graph that can be obtained from a subgraph of $G$ by contracting edges.
It is clear that $G$ contains $H$ as a topological minor implies that $G$ contains $H$ as a minor and as an immersion, but the minor containment and immersion containment are incomparable.
However, immersion, minor and topological minor are equivalent for subcubic graphs.

Immersion was considered by Nash-Williams \cite{n} when he worked on well-quasi-ordering theory.
The study of well-quasi-ordering is stimulated by a conjecture of V\'{a}zsonyi stating that subcubic graphs are well-quasi-ordered by the topological minor relation.
It is known that V\'{a}zsonyi's conjecture cannot be generalized to all graphs, as there exists an infinite antichain with respect to the topological minor relation.
In contrast, Nash-Williams conjectured that graphs are well-quasi-ordered by the immersion relation and hence provided a possible generalization of V\'{a}zsonyi's conjecture to all graphs.
Nash-Williams' conjecture was proved by Robertson and Seymour \cite{rs XXIII} by reducing it to a strengthening of their famous Graph Minor Theorem.

The cornerstone of Robertson and Seymour's proof of the Graph Minor Theorem is a structure theorem that describes, for any fixed graph $H$, the structure of $H$-minor free graphs \cite{rs XVI}.
Their theorem is stated in terms of tree-decomposition.
Roughly speaking, a tree-decomposition of a graph $G$ describes how to construct $G$ by starting with a ``piece'' of $G$ by repeatedly gluing other ``pieces'' of $G$ in a ``tree-like'' fashion.
The main result in \cite{rs XVI} states that if $G$ does not contain $H$ as a minor, then $G$ admits a tree-decomposition such that every ``piece'' can be ``nearly embedded''\footnote{The formal definition of nearly embeddable graphs is complicated. Because this formal definition is not required to understand this paper, we omit it in this paper and refer readers to \cite{rs XVI}.} in a surface in which $H$ cannot be embedded.
In fact, this statement can be considered as a global version of the decomposition theorem.
There is a local version (also proved in \cite{rs XVI}) which states that for every ``highly connected subgraph'' (or more precisely, tangle) of $G$, there is a tree-decomposition such that the underlying tree is a star and the central piece contains this ``highly connected subgraph'' and is ``nearly embeddable'' in a surface in which $H$ cannot be embedded. 
In either version, every $H$-minor free graph can be decomposed into pieces that are ``nearly simpler'' than $H$ in the sense of its Euler genus.

Decomposition theorems for excluding topological minors were developed in the same line: there are global versions \cite{d,gm} and a local version \cite{lt} in the literature.
The first such theorem was a global version proved by Grohe and Marx \cite{gm} stating that every $H$-topological minor free graph admits a tree-decomposition such that each small ``piece'' is either of ``nearly bounded maximum degree'' or ``nearly embeddable'' in a surface of bounded Euler genus.
Grohe and Marx used this to derive algorithmic results such as showing that for every graph $H$, the graph isomorphism problem on $H$-topological minor free graphs is polynomial time solvable \cite{gm}.
However, unlike Robertson and Seymour's decomposition theorem for minors, Grohe and Marx's theorem does not ensure that each piece in the decomposition is nearly simpler than $H$.
Dvo\v{r}\'{a}k \cite{d} remedies part of this drawback by proving that each piece can be made either having ``nearly bounded maximum degree'' or ``nearly embeddable'' in a surface in a way that is ``nearly impossible'' for $H$.
Inspired by Dvo\v{r}\'{a}k's work, the author and Thomas \cite{lt} proved a local version of the decomposition theorem stating that every $H$-topological minor free graph can be ``decomposed'' such that each ``piece'' either has ``maximum degree'' smaller than the maximum degree of $H$ or is ``nearly embeddable'' in a surface in a way that is ``nearly impossible'' for $H$.
Hence a result showing that each piece is ``nearly simpler'' than $H$ is obtained.
This improvement is a crucial step in resolving some open problems about topological minors, such as in a proof of Robertson's conjecture on well-quasi-ordering \cite{l_thesis,lt_wqo}, an Erd\H{o}s-P\'{o}sa type result \cite{l_half_topo} implying Thomas' conjecture, and a linear upper bound for the clustered coloring version of Haj\'{o}s' conjecture \cite{lw_topo}.

In general, each of the local version and global version has its own advantages.
The local version offers more detailed structure information and better quantitative bounds than the global version so that many applications (such as \cite{l_thesis,l_half_topo,lt_wqo,lw_minor,lw_topo,rs XX}) require the local version.
The global version sacrifices some structure information but is able to encode other important information with respect to different ``highly connected subgraphs'' into one tree-decomposition so that it is easier to be exploited in inductive arguments and provide efficient algorithms.
For example, see \cite{dfht,dhk,ddosrsv,gm}. 

A very clean local version of a decomposition theorem for excluding any fixed graph as an immersion in graphs with no edge-cut of order 3 has been obtained by the author and successfully applied to prove an Erd\H{o}s-P\'{o}sa type result for packing and covering immersions \cite{l_ep_imm}.
As the minor relation and immersion relation are equivalent for subcubic graphs, it can be shown that once edge-cuts\footnote{Edge-cuts are defined in Definition \ref{def_many}.} of order 3 are allowed, any sufficiently informative decomposition theorem with respect to immersion must be at least as complicated as Robertson and Seymour's decomposition theorem with respect to minors \cite{rs XVI}, so the complicated notion of nearly embedding is unavoidable.
On the other hand, a global version for excluding a $K_t$-immersion was proved by Wollan \cite{w_imm} stating that every $K_t$-immersion free graph can be decomposed into ``pieces'' that have a bounded number of vertices of large degree.
So as Grohe and Marx's theorem for excluding topological minors, each piece in Wollan's theorem is not ensured to be ``nearly simpler'' than $K_t$.

\subsection{Main result}

The main result of this paper is a global decomposition theorem for $H$-immersion free graphs that assembles all nice parts in the above discussion.
That is, we prove that for every graph $H$, every $H$-immersion free graph with no edge-cut of order 3 can be ``globally decomposed'' into ``pieces'' that are ``nearly simpler'' than $H$, without requiring complicated descriptions for the ``simplicity''.

First, we need the correct notion for ``decompositions''.
Global decomposition theorems for minors \cite{rs XVI} and topological minors \cite{d,gm} are stated in terms of tree-decompositions.
That is, it concerns how to break graphs by using vertex-cuts.
However, it is not effective when considering immersions.
For example, tree-decomposition is not able to distinguish the path on $t$ vertices from the graph obtained from a path on $t$ vertices by duplicating each edge $t^2$ times, where the former does not contain any graph with minimum degree at least three as an immersion, and the later contains a $K_t$-immersion\footnote{One can easily modify this example by duplicating edges and subdividing edges to obtain two simple graphs with tree-width 2, where one graph contains a $K_t$-immersion but the other graph does not contain a $K_4$-immersion.}.
Hence the decomposition theorem for immersion should address edge-cuts instead of vertex-cuts.
This leads to the notion of tree-cut decomposition which was introduced by Wollan \cite{w_imm}.

\begin{dfn}
{\rm Let $G$ be a graph.
A {\it tree-cut decomposition} of $G$ is a pair $(T,\X)$ such that $T$ is a tree, and $\X$ is a collection $(X_t: t \in V(T))$ of pairwise disjoint (not necessarily non-empty) subsets of $V(G)$ such that $\bigcup_{t \in V(T)}X_t = V(G)$.
In addition,
	\begin{itemize}
		\item for every $t \in V(T)$, the set $X_t$ is called the {\it bag} at $t$;
		\item for every subset $S$ of $V(T)$, we define $X_S$ to be $\bigcup_{t \in S}X_t$; for every subgraph $S$ of $T$, we define $X_S$ to be $\bigcup_{t \in V(S)}X_t$;
		\item for any edge $t_1t_2$ of $T$, the {\it adhesion set} of $t_1t_2$ in $(T,\X)$, denoted by $\adh_{(T,\X)}(t_1t_2)$, is the set of edges of $G$ with one end in $X_{T_1}$ and one end in $X_{T_2}$, where $T_1$ and $T_2$ are the components of $T-t_1t_2$;
		\item the {\it adhesion} of $(T,\X)$ is $\max_{e \in E(T)} \lvert \adh_{(T,\X)}(e) \rvert$.
	\end{itemize}}
\end{dfn}

Second, we need the correct notion for ``pieces''.
The ``pieces'' in the decomposition theorems for minors \cite{rs XVI} and topological minors \cite{d,gm} are the torsos of the tree-decomposition.
It can be shown that any graph is a vertex-sum of the torsos of its tree-decomposition.
That is, the graph can be constructed by repeatedly gluing together the torsos of its tree-decomposition along vertices in a certain way.
The analog of vertex-sums for tree-cut decomposition is called edge-sums. 
And the torsos of a tree-cut decomposition are defined as follows so that any graph is an edge-sum of the torsos. 

\begin{dfn}
{\rm Let $(T,\X)$ be a tree-cut decomposition of a graph $G$, where $\X=(X_t: t \in V(T))$.
For every $t$ which is a node of $T$ or a connected subgraph of $T$,
	\begin{itemize}
		\item the {\it torso} of $(T,\X)$ at $t$ is the graph obtained from $G$ by, for each edge $e$ of $T$ with exactly one end in $t$, identifying $X_{T_{t,t'}}$ into a vertex and deleting all loops incident with this new vertex, where $t'$ is the end of $e$ other than $t$, and $T_{t,t'}$ is the component of $T-t$ containing $t'$, and
		\item each vertex in the torso at $t$ but not in $X_t$ is called a {\it peripheral vertex}.
	\end{itemize}
Note that for every $t \in V(T)$ and every edge $e$ of the torso at $t$, $e$ corresponds to an edge $e'$ of $G$ such that there exists no component $T'$ of $T-t$ such that $X_{T'}$ contains all ends of $e'$.
So we may view each edge of the torso at $t$ as an edge of $G$ if there is no danger for creating confusion.
}
\end{dfn}

\begin{dfn} \label{def_edge_sum}
{\rm Let $k$ be a nonnegative integer.
Let $G_1$ and $G_2$ be graphs.
A graph $G$ is obtained from $G_1$ and $G_2$ by a {\it $k$-edge-sum} if for $i \in \{1,2\}$, $G_i$ contains a vertex $v_i$ incident with exactly $k$ non-loop edges, and there exists a bijection $f$ between the set of the $k$ non-loop edges incident with $v_1$ and the set of the $k$ non-loop edges incident with $v_2$ such that $G$ is obtained from $(G_1-v_1) \cup (G_2-v_2)$ by, for each non-loop edge $e$ incident with $v_1$, adding an edge $(e-\{v_1\}) \cup (f(e)-\{v_2\})$.
A graph is obtained from $G_1$ and $G_2$ by a {\it $(\leq k)$-edge-sum} if it is obtained from $G_1$ and $G_2$ by a $k'$-edge-sum for some nonnegative integer $k'$ with $k' \leq k$.}
\end{dfn}

It is straightforward to see that for every graph $G$ and every tree-cut decomposition of adhesion at most $k$, $G$ can be obtained from the torsos by $(\leq k)$-edge-sums.

Third, we need a measure for the ``simplicity''.
It is easy to see that if a graph $G$ contains another graph $H$ as an immersion, then for every nonnegative integer $d$, the number of vertices of $G$ with degree at least $d$ must be at least the number of vertices of $H$ with degree at least $d$.
This degree statistic provides a measure of the simplicity.
Theorem \ref{global_intro} is the main theorem of this paper, stating that the torsos nearly violate the above necessary degree condition and hence are ``nearly simpler'' than $H$.
For a technical reason, we require that the graph $H$ is not ``too simple''.
See Theorem \ref{global} for a more comprehensive version of Theorem \ref{global_intro}.

A graph is {\it exceptional} if it contains exactly one vertex of degree at least two, and this vertex is incident with a loop.
Note that subdividing any edge of an exceptional graph makes it non-exceptional, and this operation enables us to apply Theorem \ref{global_intro} in many applications.

\begin{theorem} \label{global_intro}
For any positive integers $d$ and $h$, there exist integers $\eta$ and $\xi$ such that the following holds.
Let $H$ be a non-exceptional graph on $h$ vertices with maximum degree $d$.
Let $G$ be a graph with no edge-cut of order exactly 3 such that $G$ does not contain an $H$-immersion.
Then there exists a tree-cut decomposition $(T,\X)$ of $G$ of adhesion at most $\eta$ such that for every $t \in V(T)$, there exists $Z_t \subseteq E(G)$ with $\lvert Z_t \rvert \leq \xi$ such that if $G_t$ is the torso at $t$, then there exists a nonnegative integer $k_t$ such that the number of vertices of degree at least $k_t$ in $G_t-Z_t$ is less than the number of vertices of degree at least $k_t$ in $H$.
\end{theorem}

Recall that if edge-cuts of order 3 are allowed, then ``nearly embedding'' must be included in the conclusion of a structure theorem for excluding immersions.
So forbidding edge-cuts of order 3 is essential in order to obtain a clean structure theorem.

\subsection{Applications}

Theorem \ref{global_intro} is expected to be useful to reduce problems on $H$-immersion free graphs to graphs whose degree sequence witnesses the $H$-immersion freeness.
As a supportive evidence, in an accompanied paper \cite{l_imm_clu}, we will use Theorem \ref{global_intro} to prove, for every graph $H$, an upper bound and a lower bound for the minimum $k$ such that every $H$-immersion free graph can be partitioned into $k$ induced subgraphs with bounded component size, where the upper and lower bounds coincide for infinitely many graphs $H$ and differ by 1 for the rest of graphs $H$.
This is the clustered coloring version of a question proposed independently by Lescure and Meyniel \cite{lm} and Abu-Khzam and Langston \cite{al} regarding the immersion-analog of Hadwiger's conjecture and Haj\'{o}s' conjecture on coloring.
Note that the clustered coloring version of Hadwiger's conjecture and Haj\'{o}s' conjecture have been extensively studied.
The result in \cite{l_imm_clu} shows that the answer for the immersion version of the clustered coloring problem behaves significantly differently from the minor and topological minor version.
We refer readers to \cite{l_imm_clu} for more details.

A simple corollary of Theorem \ref{global_intro} is the following statement whose $d$-edge-connected case is a result of Marx and Wollan \cite[Theorem 1.2]{mw}. 

\begin{corollary} \label{mw_thm_intro}
For any positive integers $d,h$ with $d \geq 4$, there exists a positive integer $\xi=\xi(d,h)$ such that the following hold.
If $G$ is an $H$-immersion free graph for some graph $H$ on at most $h$ vertices with maximum degree at most $d$, and every component of $G$ has no edge-cut of order $k$ for every $3 \leq k \leq d-1$, then there exists a tree-cut decomposition $(T,\X)$ of $G$ of adhesion at most $\xi$ such that 
	\begin{enumerate}
		\item the maximum degree of $T-S$ is at most $\xi$, where $S = \{e \in E(T): \lvert \adh_{(T,\X)}(e) \rvert \leq 2\}$, and
		\item $\lvert X_t \rvert \leq \xi$ for every $t \in V(T)$.
	\end{enumerate}
In particular, the tree-cut width\footnote{The notion of tree-cut width is defined in Definition \ref{def_tree_cut_width}.} of $G$ is at most $2\xi$. 
\end{corollary}

Corollary \ref{mw_thm_intro} can be used to prove the following enumerative result. 

\begin{theorem} \label{counting_basic_intro}
For any positive integers $d,h$ with $d \geq 4$, there exists a positive integer $c$ such that the following hold.
Let $H$ be a graph on at most $h$ vertices with maximum degree at most $d$.
Let $\C$ be the set of (vertex-)unlabelled $H$-immersion free graphs whose every maximal 2-edge-connected subgraph is $d$-edge-connected.
Then 
	\begin{enumerate}
		\item for every positive integer $m$, the number of members of $\C$ with $m$ edges and with no isolated vertex is at most $c^m$, and
		\item for every positive integer $n$, the number of simple members of $\C$ with $n$ vertices is at most $c^n$.
	\end{enumerate}
\end{theorem}

Note that the $d$-edge-connectivity in Theorem \ref{counting_basic_intro} is required as there are superexponentially many unlabelled $(d-1)$-regular $(d-1)$-edge-connected simple graphs on $n$ vertices \cite{bc,l_random,w_counting} (see \cite{cls}), and every $(d-1)$-regular graph does not contain $H$ as an immersion if $H$ has maximum degree $d$.
Analogous results for minors and topological minors are known in the literature: for every graph $H$, there exists $c$ such that there are at most $c^n$ unlabelled simple $n$-vertex $H$-minor free graphs \cite{afs,dn}, implying that there are at most $n!c^n$ labelled simple $n$-vertex $H$-minor free graphs \cite{nstw}; for any integers $d,d'$ and graph $H$ with maximum degree at most $d$, there exists $c$ such that there are at most $c^n$ unlabelled simple $d$-connected $n$-vertex $H$-topological minor free graphs with maximum degree at most $d'$ \cite{cls}.

Such enumeration results are motivated by the work in \cite{msw} about small addable classes.
A class $\C$ of graphs is {\it small} if there exists a constant $c$ such that for every integer $n$, there are at most $n!c^n$ labelled simple $n$-vertex graphs in $\C$.
So Theorem \ref{counting_basic_intro} and the aforementioned results in \cite{afs,cls,dn,nstw} show that certain classes are small.
A class $\C$ of graphs is {\it addable} if 
	\begin{itemize}
		\item $G \in \C$ if and only if every component of $G$ is in $\C$, and
		\item if $G_1,G_2 \in \C$, then the graph obtained from a disjoint union of $G_1$ and $G_2$ by adding an edge between $G_1$ and $G_2$ is in $\C$.
	\end{itemize}
McDiarmid, Steger and Welsh \cite{msw} proved many results about small addable classes of simple graphs, such as the following.

\begin{theorem}[{\cite{msw}}] \label{small_addable}
Let $\C$ be a small addable class of simple graphs.
Then the following hold.
	\begin{enumerate}
		\item $\lim_{n \to \infty}(N(n)/n!)^{1/n} = c$ for some constant $c$, where $N(n)$ is the number of $n$-vertex labelled graphs in $\C$.
		\item For every positive integer $k$, if $K_{1,k+1} \in \C$, then there exist constants $b$ and $n_0$ such that for every $n \geq n_0$, selecting a graph $G$ uniformly from  the $n$-vertex graphs in $\C$, the probability that $G$ has fewer than $a_kn$ vertices of degree $k$ is at most $e^{-a_kn}$, where $a_k = b/(c^k(k+2)!)$.
		\item The probability that the graph $G$ in Statement 2 has an isolated vertex is at least $a_1/e+o(1)$.
	\end{enumerate}
\end{theorem}

As Theorem \ref{counting_basic_intro} shows that the class $\C$ mentioned in Theorem \ref{counting_basic_intro} is small and addable, all conclusions of Theorem \ref{small_addable} apply to simple graphs in $\C$.

Now we briefly discuss the proof of Theorem \ref{global_intro}.
Though the global decomposition theorem for excluding minors can be easily derived from the local version \cite{rs XVI}, it is unclear how to use similar arguments to derive Theorem \ref{global_intro} from the local decomposition in \cite{l_ep_imm}.
So we use a strategy different from \cite{rs XVI} to derive Theorem \ref{global_intro} from the results in \cite{l_ep_imm}.
A more detailed proof sketch is included in Section \ref{sec:sketch}.

A byproduct of our proof is a duality theorem for maximum order of edge-tangles and tree-cut torso-width.
The notion of edge-tangles will be defined in Section \ref{sec:smooth}.

\begin{dfn}
{\rm Let $G$ be a graph.
	\begin{itemize}
		\item The {\it torso-width} of a tree-cut decomposition $(T,\X)$ of $G$ is the minimum $w$ such that for every $t \in V(T)$, the torso at $t$ in $(T,\X)$ has at most $w$ edges.
		\item The {\it tree-cut torso-width} of $G$ is the minimum $w$ such that $G$ admits a tree-cut decomposition of torso-width $w$.
	\end{itemize}}
\end{dfn}

The notion of tree-cut torso-width is natural, but it seems that it was not considered in the literature in our knowledge.
Tree-cut torso-width is closely related to carving-width which was introduced by Seymour and Thomas \cite{st}.
(See Section \ref{sec:app} for a formal definition of carving-width.)
Carving width can be viewed as an edge-analog of branch width which is another extensively studied width parameter. 
Robertson and Seymour \cite{rs X} proved that having bounded branch width is equivalent to having bounded tree-width and is equivalent to having no tangle of large order.

Tree-cut torso-width can be viewed as an edge-analog of tree-width based on their definitions.
So one might expect that having bounded carving-width is equivalent to having bounded torso-width.
It is easy to show that it is indeed the case if graphs are loopless; however, it is not true when loops are allowed, as adding loops does not change the carving-width but can make tree-cut torso-width arbitrarily large.
Same situation happens for the relationship between edge-tangles, carving-width and tree-cut torso-width.
One corollary of general work of Diestel and Oum on abstract separation systems \cite{do2,do} shows that for loopless graphs, having no edge-tangle of large order is equivalent to having bounded carving-width.
A byproduct of our proof of Theorem \ref{global_intro} is Proposition \ref{width_relation_intro} which shows that it is true even when loops are allowed and gives an independent proof for the case of loopless graphs. 

\begin{proposition} \label{width_relation_intro}
Let $G$ be a graph.
Let $w$ be a nonnegative integer.
	\begin{enumerate}
		\item If the tree-cut torso-width of $G$ is at most $w$, then the carving-width of $G$ is at most $w$ and there exists no edge-tangle of order $w+1$ in $G$.
		\item If the carving-width of $G$ is at most $w$ and $G$ is loopless, then the tree-cut torso-width of $G$ is at most $3w/2$.
		\item If there exists no edge-tangle of order $w$ in $G$, then the tree-cut torso-width is at most $3w-3$.
	\end{enumerate}
\end{proposition}

We conclude this subsection by discussing the relationship between our results and tree-cut width which is the main width parameter for tree-cut decomposition considered in the literature introduced by Wollan \cite{w_imm}.

\begin{dfn} \label{def_tree_cut_width}
{\rm The {\it width} of a tree-cut decomposition $(T,\X)$ of a graph $G$ is the maximum among the adhesion of $(T,\X)$ and $\min_{t \in V(T)}\lvert V(\bar{H}_t) \rvert$, where $\bar{H}_t$ is the graph, called the {\it 3-center at $t$}, obtained from the torso at $t$ by repeatedly suppressing peripheral vertices of degree at most two until every peripheral vertex has degree at least 3.
The {\it tree-cut width} of $G$ is the minimum width of a tree-cut decomposition of $G$.}
\end{dfn}
The motivation for considering 3-centers is to avoid very simple graphs such as stars having large width.
Tree-cut width is an effective width parameter with respect to problems about immersions.
For example, Wollan \cite{w_imm} proved an analogy of the Grid Minor Theorem: every graph with large tree-cut width contains a large wall as an immersion.

It is easy to see that every graph with bounded tree-cut torso-width has bounded tree-cut width, but the converse is not true.
Hence large tree-cut torso-width does not ensure the existence of a large wall immersion.
But this issue can be easily fixed.
Note that for any fixed integer $k$, if $G$ is a 3-edge-connected graph such that every vertex is incident with at most $k$ loops and every pair of vertices has at most $k$ parallel edges between them, then $G$ has bounded tree-cut width implies that $G$ has bounded tree-cut torso-width.
This extra assumption can be made when considering immersion problems.
As for any fixed graph $H$, if there exists an $H$-immersion in a graph $G$, then there exists an $H$-immersion in $G$ such that only at most $\lvert E(H) \rvert$ parallel edges between each fixed pair of vertices of $G$ and at most $\lvert E(H) \rvert$ loops of $G$ are involved in the $H$-immersion.
So to test whether a graph contains an $H$-immersion or not, one can only keep at most $\lvert E(H) \rvert$ loops incident with each vertex and at most $\lvert E(H) \rvert$ parallel edges between each pair of vertices.
In addition, to test if a graph $G$ contains a 3-edge-connected graph $H$ or a wall as an immersion, it suffices to first decompose $G$ along edge-cuts of order at most 2 into smaller graphs with the edges in the edge-cuts split off, and then test whether each smaller graph contains an $H$-immersion or not.
Hence Wollan's wall-immersion theorem can be adapted to tree-cut torso-width of the modified graphs as well.

\subsection{Proof sketch and organization} \label{sec:sketch}

In this subsection, we sketch the proof of Theorem \ref{global_intro}.
Other results (Corollary \ref{mw_thm_intro}, Theorem \ref{counting_basic_intro}, and Proposition \ref{width_relation_intro}) are fairly simple corollaries of Theorem \ref{global_intro} or lemmas developed on the way for proving it and are formally proved in Section \ref{sec:app}.

Let $G$ and $H$ be the graphs stated in Theorem \ref{global_intro}.
Note that given a tree-cut decomposition $(T,\X)$ of $G$, for each edge $e$ of $T$, it defines an edge-cut of $G$, where each side in the edge-cut is the union of the bags of the nodes in a component of $T-e$.
And for each node $t$ of $T$, the edge-cuts of $G$ defined by the edges of $T$ incident with $t$ defines a cross-free family if we always put the side having the bag at $t$ into the second part of each edge-cut.
(Cross-free families will be formally defined in Section \ref{sec:cross-free}.)
Note that the torso at $t$ can be easily told by this cross-free family.

To show that a tree-cut decomposition $(T,\X)$ satisfies Theorem \ref{global_intro}, our attention is on the nodes whose torsos have sufficiently many edges, as for every node $t$ whose torso has only few edges, we can put all edges of the torso into $Z_t$.
Hence we would like to show that each node whose torso has many edges corresponds to an edge-tangle of large order and show that this torso comes from the cross-free family obtained by the local decomposition theorem with respect to this edge-tangle.

We first discuss the correspondence between nodes and edge-tangles.
Edge-tangles is formally defined in Section \ref{sec:smooth}.
Roughly speaking, each edge-tangle comes from a choice of a side of each edge-cut $[A,B]$ by indicating which $A$ or $B$ is ``more important''.
For a node $t$ of $T$ whose torso has many edges, if for each edge-cut $[A,B]$ of order less than a constant $\theta$, one side of $[A,B]$ contains only few edges of the torso, then we can say the other side is more important, so we can define an edge-tangle $\E$.
In this case, we can say that the node $t$ corresponds to $\E$.
We call a tree-cut decomposition with this property a $\theta$-smooth tree-cut decomposition.
See Section \ref{sec:smooth} for a formal definition.
The existence of a $\theta$-smooth tree-cut decomposition is proved in Theorem \ref{theta-smooth}.
The formal correspondence between edge-tangles and nodes whose torsos have many edges is also proved in Section \ref{sec:smooth}.

Note that it is easy to see how an edge-cut given by an edge of $T$ distinguishes different edge-tangles given by nodes of $T$.
But we need to understand other edge-cuts that distinguish different edge-tangles for the future usage.
This is the motivation of ``separators'' defined in Section \ref{sec:smooth}.
Partial information about such separators can also be told from a $\theta$-smooth tree-cut decomposition, as shown in Section \ref{sec:smooth}.

Now we know that each node whose torso with many edges corresponds to an edge-tangle.
So for each edge-tangle $\E_t$ defined by such a node $t$ of $T$, we can apply the local decomposition theorem (Lemma \ref{basic_exclu_imm}) with respect to $\E_t$ to obtain a cross-free family $\D_t$ such that if the collection of the edge-cuts given by the edges of $T$ incident with $t$ is $\D_t$, then the torso at $t$ would satisfy Theorem \ref{global_intro}.
However, edge-cuts in $\D_{t_1}$ can cross edge-cuts in $\D_{t_2}$ for distinct nodes $t_1$ and $t_2$, so it is unlikely that we can ``realize'' $\D_t$ for all nodes $t$ simultaneously in $(T,\X)$.
So we would like to modify each $\D_t$ to obtain a new cross-free family $\D_t'$ such that the edge-cuts in $\D_{t_1}'$ do not cross the edge-cuts in $\D_{t_2}'$ for different nodes $t_1$ and $t_2$, and Theorem \ref{global_intro} would still hold if we can realize $\D_t'$ in $(T,\X)$ instead of $\D_t$.

To make sure that the edge-cuts in $\D_{t_1}$ do not cross the edge-cuts in $\D_{t_2}$ for distinct nodes $t_1,t_2$, we need to consider edge-cuts that distinguish the edge-tangles consistent with $\D_{t_1}$ from the edge-tangles consistent with $\D_{t_2}$.
This is the motivation of ``segregators'' and ``guards'' defined in Section \ref{sec:cross-free}.
In Section \ref{sec:cross-free}, we develop tools to show that we can indeed modify each $\D_t$ into $\D_t'$ such that $\D_t'$ still preserves certain nice properties of $\D_t$ such that realizing $\D'_t$ in a tree-cut decomposition still ensures the validity of Theorem \ref{global_intro}, and for distinct nodes $t_1$ and $t_2$, the edge-cuts in $\D_{t_1}'$ do not cross the edge-cuts in $\D_{t_2}'$.

In Section \ref{sec:4-edge-conn}, we formally show how to modify a smooth tree-cut decomposition $(T,\X)$ to construct a new tree-cut decomposition $(T^*,\X^*)$ so that each node $t$ of $T^*$ whose torso in $(T^*,\X^*)$ has many edges comes from a node of $T$ whose torso in $(T,\X)$ has many edges, and $\D_t'$ is realized in $(T^*,\X^*)$ to complete the proof.

\subsection{Definitions} \label{sec:def}

The following notions and notations will be frequently used in this paper.

\begin{dfn} \label{def_many}
{\rm Let $G$ be a graph.
	\begin{itemize}
		\item An {\it edge-cut} $[A,B]$ of a graph $G$ is an ordered pair of disjoint subsets of $V(G)$ such that $A \cup B = V(G)$.
The {\it order} of an edge-cut $[A,B]$ of $G$ is the number of edges of $G$ with one end in $A$ and one end in $B$.
		\item Let $(T,\X)$ be a tree-cut decomposition of $G$.
Let $t$ be a node of $T$ or a connected subgraph of $T$.
Let $e$ be an edge of $T$ with at most one end in $t$.
We define $[A_{e,t},B_{e,t}]$ to be the edge-cut of $G$ with $B_{e,t} = \bigcup_{t''} X_{t''}$, where the union is over all nodes $t''$ contained in the component of $T-e$ containing $t$.
		\item For every subset $S$ of $V(G)$, we define $G[S]$ to be the subgraph of $G$ induced by $S$.
		\item The {\it degree sequence} of $G$ is the non-increasing sequence $(d_1,d_2,...,d_{\lvert V(G) \rvert})$ such that there exists a bijection $\iota: V(G) \rightarrow [\lvert V(G) \rvert]$ such that for every $v \in V(G)$, the degree of $v$ equals $d_{\iota(v)}$.
	\end{itemize}
	}
\end{dfn}

\section{Smooth tree-cut decompositions} \label{sec:smooth}

Let $(T,\X)$ be a tree-cut decomposition of a graph $G$.
For every positive integer $k$, a {\it pseudo-$k$-cell} in $(T,\X)$ is a component $C$ of the forest obtained from $T$ by deleting all edges of $T$ whose adhesion set has size less than $k$; a {\it $k$-cell} in $(T,\X)$ is a pseudo-$k$-cell $C$ such that the number of edges in the torso at $C$ is at least $k$.

Let $\theta$ be a positive integer.
Let $(T,\X)$ be a tree-cut decomposition of a graph $G$.
We say that $(T,\X)$ is {\it $\theta$-smooth} if for every $C$ which is a $\theta$-cell in $(T,\X)$ or a node of $T$, for any sets $Y$ and $Z$ of edges of the torso at $C$ with $\lvert Y \rvert = \lvert Z \rvert \leq \theta$, there does not exist an edge-cut $[A,B]$ of $G$ of order less than $\lvert Y \rvert$ such that every edge in $Y$ is incident with $A$ and every edge in $Z$ is incident with $B$. 

The following is an easy but useful restatement of the definition of the $\theta$-smooth property. 

\begin{lemma} \label{equivdefsmooth}
Let $(T,\X)$ be a tree-cut decomposition of a graph $G$.
Let $\theta$ be a positive integer.
Then $(T,\X)$ is $\theta$-smooth if and only if for every $C$ that is a node of $T$ or a $\theta$-cell in $(T,\X)$, there exists no edge-cut $[A,B]$ of $G$ of order less than $\theta$ such that each of $A$ and $B$ is incident with at least $\lvert [A,B] \rvert+1$ edges of the torso at $C$.
\end{lemma}

\begin{pf}
	Assume that there exists $C$ that is a node of $T$ or a $\theta$-cell in $(T,\X)$ and there exists an edge-cut $[A,B]$ of $G$ of order less than $\theta$ such that each $A$ and $B$ is incident with at least $\lvert [A,B] \rvert+1$ edges of the torso at $C$.
	Then there exists a set $Y$ of edges of the torso at $C$ incident with $A$ with $\lvert Y \rvert =\lvert [A,B] \rvert+1 \leq \theta$, and there exists a set $Z$ of the edges of the torso at $C$ incident with $B$ with $\lvert Z \rvert = \lvert [A,B] \rvert+1$.
	Note that $\lvert [A,B] \rvert < \lvert Y \rvert$.
	Hence $(T,\X)$ is not $\theta$-smooth.
	
	Assume that $(T,\X)$ is not $\theta$-smooth.
	Then there exists $C$ which is a $\theta$-cell in $(T,\X)$ or a node of $T$, and there exist sets $Y$ and $Z$ of edges of the torso at $C$ with $\lvert Y \rvert = \lvert Z \rvert \leq \theta$ and an edge-cut $[A,B]$ of $G$ of order less than $\lvert Y \rvert$ such that every edge in $Y$ is incident with $A$ and every edge in $Z$ is incident with $B$.
	So $[A,B]$ is an edge-cut of $G$ of order less than $\lvert Y \rvert \leq \theta$ such that $A$ is incident with at least $\lvert Y \rvert \geq \lvert [A,B] \rvert+1$ edges and $B$ is incident with at least $\lvert Z \rvert \geq \lvert [A,B] \rvert+1$ edges.
\end{pf}

\bigskip

Let $(T,\X)$ be a tree-cut decomposition of a graph $G$.
For integers $i,j$ with $1 \leq i \leq \lvert E(G) \rvert$ and $1 \leq j \leq \lvert E(G) \rvert$, we define 
	\begin{itemize}
		\item $a_{i,j}$ to be the number of $i$-cells $L$ of $(T,\X)$ such that the torso at $L$ has at least $j$ edges, and
		\item $a_i$ is the sequence $(a_{i,\lvert E(G) \rvert},a_{i,\lvert E(G) \rvert-1},...,a_{i,1})$.
	\end{itemize}
Let $k$ be a positive integer.
The {\it $k$-signature} of $(T,\X)$ is the sequence $(a_k,a_{k-1},...,a_1)$. 

Intuitively, if a tree-cut decomposition of $G$ is not $\theta$-smooth, then we can ``insert'' an edge-cut of $G$ into this tree-cut decomposition to break a cell $L$ into two cells with fewer edges, so it decreases some entry in the signature with the price that increases entries appearing later than the aforementioned entry in the signature; and if one can make sure that every entry appearing earlier than the aforementioned entry does not increase, then the lexicographic order of the signature decreases.
The following theorem shows that this intuition is correct, if we choose the edge-cut that breaks the cell carefully.

\begin{theorem} \label{theta-smooth}
Let $G$ be a graph.
Let $\theta$ be a positive integer.
Let $(T,\X)$ be a tree-cut decomposition of $G$ of lexicographically minimum $\theta$-signature.
Then $(T,\X)$ is $\theta$-smooth.
\end{theorem}

\begin{pf}
Suppose to the contrary that $(T,\X)$ is not $\theta$-smooth.
By Lemma \ref{equivdefsmooth}, there exist $L$ that is either a $\theta$-cell in $(T,\X)$ or a node of $T$, a positive integer $k \in [\theta]$, and an edge-cut $[A,B]$ of $G$ of order $k-1$ such that there are at least $\lvert [A,B] \rvert+1=k$ edges in the torso at $L$ incident with $A$ and there are at least $\lvert [A,B] \rvert+1=k$ edges in the torso of $L$ incident with $B$.
We further choose $L$ and $[A,B]$ such that $\lvert \{e \in E(T): e$ has at most one end in $V(L)$, either $A_{e,L} \subseteq A$ or $A_{e,L} \subseteq B\} \rvert$ is as large as possible.

\noindent{\bf Claim 1:} For every edge $e$ of $T$ with at most one end in $V(L)$, either
		\begin{itemize}
			\item $A_{e,L} \subseteq A$ or $A_{e,L} \subseteq B$, or
			\item $\lvert [A \cap A_{e,L}, B \cup B_{e,L}] \rvert < \lvert [A_{e,L},B_{e,L}] \rvert$ and $\lvert [A \cup B_{e,L}, B \cap A_{e,L}] \rvert < \lvert [A_{e,L},B_{e,L}] \rvert$.
		\end{itemize}
		
\noindent{\bf Proof of Claim 1:}
Suppose to the contrary that there exists an edge $e$ of $T$ with at most one end in $V(L)$ such that $A_{e,L} \not \subseteq A$, $A_{e,L} \not \subseteq B$, and either $\lvert [A \cap A_{e,L}, B \cup B_{e,L}] \rvert \geq \lvert [A_{e,L},B_{e,L}] \rvert$ or $\lvert [A \cup B_{e,L}, B \cap A_{e,L}] \rvert \geq \lvert [A_{e,L},B_{e,L}] \rvert$.
By symmetry, we may assume that $\lvert [A \cap A_{e,L}, B \cup B_{e,L}] \rvert \geq \lvert [A_{e,L},B_{e,L}] \rvert$.
	
Let $S_1$ be the set of the edges of the torso at $L$ incident with $B$ but not incident with $B \cap B_{e,L}$.
Since $e$ has at most one end in $V(L)$, no edge of the torso at $L$ has both ends in $A_{e,L}$.
So every edge in $S_1$ has one end in $B \cap A_{e,L}$ and one end in $A \cap B_{e,L}$.
Let $S_2$ be the set of edges of $G$ with one end in $B \cap A_{e,L}$ and one end in $A \cap B_{e,L}$.
So $S_1 \subseteq S_2$.
Since there are at least $\lvert [A,B] \rvert+1$ edges of the torso at $L$ incident with $B$, there are at least $\lvert [A,B] \rvert+1 - \lvert S_1 \rvert$ edges of the torso at $L$ incident with $B \cap B_{e,L}$.
	
In addition, $\lvert [A,B] \rvert + \lvert [A_{e,L},B_{e,L}] \rvert = \lvert [A \cap A_{e,L}, B \cup B_{e,L}] \rvert + \lvert [A \cup A_{e,L}, B \cap B_{e,L}] \rvert + 2\lvert S_2 \rvert \geq \lvert [A_{e,L}, B_{e,L}] \rvert + \lvert [A \cup A_{e,L}, B \cap B_{e,L}] \rvert + 2\lvert S_2 \rvert$.
So $\lvert [A \cup A_{e,L}, B \cap B_{e,L}] \rvert \leq \lvert [A,B] \rvert - 2\lvert S_2 \rvert \leq \lvert [A,B] \rvert - \lvert S_1 \rvert \leq \lvert [A,B] \rvert$.
Since $A \cup A_{e,L} \supseteq A$, there are at least $\lvert [A,B] \rvert+1 \geq \lvert [A \cup A_{e,L}, B \cap B_{e,L}] \rvert+1$ edges of the torso at $L$ incident with $A \cup A_{e,L}$.
And there are at least $\lvert [A,B] \rvert+1 - \lvert S_1 \rvert \geq \lvert [A \cup A_{e,L}, B \cap B_{e,L}] \rvert+1$ edges of the torso at $L$ incident with $B \cap B_{e,L}$.

Let $S = \{e' \in E(T): e'$ has at most one end in $V(L)$, either $A_{e',L} \subseteq A$ or $A_{e',L} \subseteq B\}$.
Since $A_{e,L} \not \subseteq A$ and $A_{e,L} \not \subseteq B$, for every edge $e' \in S$, $A_{e,L} \not \subseteq A_{e',L}$.
So for every $e' \in S$, either $A_{e',L} \subseteq A_{e,L}$ or $A_{e',L} \subseteq B_{e,L}$, so either $A_{e',L} \subseteq A \cup A_{e,L}$ or $A_{e',L} \subseteq B \cap B_{e,L}$.
In addition, $A_{e,L} \subseteq A \cup A_{e,L}$ and $e \not \in S$.
Therefore, $[A \cup A_{e,L},B \cap B_{e,L}]$ is a better choice than $[A,B]$, a contradiction. 
$\Box$

Since there are at least $\lvert [A,B] \rvert+1 =k$ edges of the torso at $L$ incident with $A$, $L$ is contained in an $r$-cell for some $r \geq k$.
Let $k^*$ be the largest integer with $k^* \in [\theta]$ such that $L$ is contained in a $k^*$-cell $L^*$.
Note that $k^* \geq k$.

If $L$ is not a $\theta$-cell, then $L$ is a node, and we define $T_0=T$ and let $t_0=L$; if $L$ is a $\theta$-cell, define $T_0$ to be the tree obtained from $T$ by contracting $L$ into a new vertex $t_0$.
Let $T^*$ be the tree obtained from a union of two disjoint copies $T',T''$ of $T_0$ by adding an edge $t_0't_0''$, where $t_0',t_0''$ are the copies of $t_0$ in $T',T''$, respectively.
For any $x$ which is a node, an edge or a subgraph of $T_0$, we denote the copy of $x$ in $T'$ and $T''$ by $x'$ and $x''$, respectively.
	
For every node $t \in V(T_0)-\{t_0\}$, we define $X^*_{t'}=X_t \cap A$ and $X^*_{t''}=X_t \cap B$. 
Define $X^*_{t_0'}=X_{L^*} \cap A$ and $X^*_{t_0''}=X_{L^*} \cap B$.
Then $(T^*,\X^*)$ is a tree-cut decomposition of $G$, where $\X^*=(X^*_t: t \in V(T^*))$.
	
For every edge $e \in E(T^*)-\{t_0't_0''\}$, let $f(e)$ be the edge of $T_0$ such that $e=f(e)'$ or $e=f(e)''$.
Note that for every $e \in E(T')$, $\adh_{(T^*,\X^*)}(e)$ is the set of edges with one end in $A \cap A_{f(e),L}$ and one end in $B \cup B_{f(e),L}$; for every edge $e \in E(T'')$, $\adh_{(T^*,\X^*)}(e)$ is the set of edges with one end in $B \cap A_{f(e),L}$ and one end in $A \cup B_{f(e),L}$.
So by Claim 1, for every edge $e \in E(T^*)-\{t_0't_0''\}$, $\lvert \adh_{(T^*,\X^*)}(e) \rvert \leq \lvert \adh_{(T,\X)}(f(e)) \rvert$.
And $\adh_{(T^*,\X^*)}(t_0't_0'')$ is the set of edges with one end in $A$ and one end in $B$, so $\lvert \adh_{(T^*,\X^*)}(t_0't_0'') \rvert = k-1<k^*$.
So for every integer $r$ with $\theta \geq r \geq k^*$ and every $r$-cell $R$ in $(T^*,\X^*)$, there exists a pseudo-$r$-cell $f(R)$ in $(T,\X)$ such that either $R \subseteq f(R)'$ or $R \subseteq f(R)''$. 
		
\noindent{\bf Claim 2:} For every integer $r$ with $\theta \geq r \geq k^*$ and every $r$-cell $R$ in $(T^*,\X^*)$ with $\{t_0',t_0''\} \cap V(R)=\emptyset$, the number of edges of the torso at $f(R)$ is at least the number of edges of the torso at $R$.
Furthermore, if the number of edges of the torso at $R$ equals the number of edges of the torso at $f(R)$, then $R$ is the only $r$-cell in $(T^*,\X^*)$ contained in $f(R)' \cup f(R)''$, and either
	\begin{itemize}
		\item $X_{f(R)} \subseteq A$ and $R=f(R)'$, or 
		\item $X_{f(R)} \subseteq B$ and $R=f(R)''$.
	\end{itemize} 
In particular, $f(R)$ is an $r$-cell in $(T,\X)$.
	
\noindent{\bf Proof of Claim 2:}
We may assume that the number of edges of the torso at $f(R)$ is at most the number of edges of the torso at $R$, for otherwise we are done.
By symmetry, we may assume that $R \subseteq T'$.
Let $w$ be the peripheral vertex of the torso at $R$ in $(T^*,\X^*)$ such that $w$ is obtained from contracting the component of $T^*-V(R)$ containing $t_0'$.
Let $e_w$ be the edge of $T^*$ with one end in $V(R)$ and one end in the component of $T^*-V(R)$ containing $t_0'$.

We first assume that $A_{f(e_w),L} \not \subseteq A$ and $A_{f(e_w),L} \not \subseteq B$.
By Claim 1, $\lvert [A \cap A_{f(e_w),L}, B \cup B_{f(e_w),L}] \rvert < \lvert [A_{f(e_w),L},B_{f(e_w),L}] \rvert = \lvert \adh_{(T,\X)}(f(e_w)) \rvert$.

Let $S_w$ be the set of edges of the torso at $R$ incident with $w$.
That is, $S_w$ consists of the edges of $G$ between $A \cap A_{f(e_w),L}$ and $B \cup B_{f(e_w),L}$.
So $\lvert S_w \rvert < \lvert \adh_{(T,\X)}(f(e_w)) \rvert$.
Note that every edge in $\adh_{(T,\X)}(f(e_w))$ is an edge of the torso at $f(R)$.
	
Let $Q_w$ be the set of edges of the torso at $R$ not incident with $w$.
Then all ends of any edge in $Q_w$ are contained in $A \cap A_{f(e_w),L}$.
So every edge in $Q_w$ is an edge of the torso at $f(R)$ and is not in $\adh_{(T,\X)}(f(e_w))$.
Hence the set of the edges of the torso at $f(R)$ contains the union of the two disjoint sets $\adh_{(T,\X)}(f(e_w))$ and $Q_w$.
So the number of edges of the torso at $f(R)$ is at least $\lvert \adh_{(T,\X)}(f(e_w)) \rvert + \lvert Q_w \rvert > \lvert S_w \rvert + \lvert Q_w \rvert$.
Therefore, the number of edges of the torso at $R$ equals $\lvert S_w \rvert + \lvert Q_w \rvert$ which is smaller than the number of edges at the torso at $f(R)$.
So we are done.

Hence we may assume that $A_{f(e_w),L} \subseteq A$ or $A_{f(e_w),L} \subseteq B$.
If $A_{f(e_w),L} \subseteq B$, then the torso at $R$ has no edge, a contradiction.
So $A_{f(e_w),L} \subseteq A$.
Hence $\adh_{(T,\X)}(f(e_w)) \subseteq \adh_{(T^*,\X^*)}(e_w)$.
Since $e_w$ has exactly one end in $V(R)$, $\lvert \adh_{(T^*,\X^*)}(e_w) \rvert<r$.
So $f(e_w)$ has exactly one end in $f(R)$.
That is, $f(e_w)$ is the edge incident with $f(R)$ such that $f(R)$ and $t_0$ belong to different components in $T_0-f(e_w)$.
Therefore, $X_{f(R)} \subseteq A_{f(e_w),L} \subseteq A$.
This implies that for every edge $e'$ of $T'$ in which $f(e')$ is an edge of $f(R)$, $\adh_{(T,\X)}(f(e')) \subseteq \adh_{(T^*,\X^*)}(e')$, so $e'$ is contained in $R$.
Therefore, $R=f(R)'$.
Similarly, since $A_{f(e_w),L} \subseteq A$, we know $A_{f(e_w)'',t_0''}=\emptyset$, so there exists no $r$-cell in $(T^*,\X^*)$ contained in $f(R)''$.

Since $R$ is an $r$-cell, the torso at $R$ contains at least $r$ edges.
Hence the torso at the pseudo-$r$-cell $f(R)$ contains at least $r$ edges.
So $f(R)$ is an $r$-cell.
$\Box$

\noindent{\bf Claim 3:} For every $r \in [\theta]$ with $r \geq k^*$, if there exists an $r$-cell $Q$ in $(T^*,\X^*)$ containing $t_0'$ or $t_0''$, then $r=k^*$, the number of edges of the torso at $Q$ in $(T^*,\X^*)$ is strictly smaller than the number of edges of the torso at $L^*$ in $(T,\X)$, and $f(Q)=L^*$. 

\noindent{\bf Proof of Claim 3:}
Let $Q$ be an $r$-cell in $(T^*,\X^*)$ for some $\theta \geq r \geq k^*$ containing $t_0'$ or $t_0''$.
Since $\lvert \adh_{(T^*,\X^*)}(t_0't_0'') \rvert=k-1<k^*$, $Q$ contains exactly one of $t_0'$ and $t_0''$.
By symmetry, we may assume that $Q$ contains $t_0'$.
So $Q$ is contained in $T'$.

We first assume that $L$ is a $\theta$-cell.
Then $k^*=\theta$ and $L=L^*$ by the definition of $k^*$ and $L^*$.
Since $k^* \leq r \leq \theta$, $r=k^*=\theta$.
By Claim 1, for every $e \in E(T^*)$ incident with $t_0'$, $\lvert \adh_{(T^*,\X^*)}(e) \rvert \leq \lvert \adh_{(T,\X)}(f(e)) \rvert <\theta$, so $Q$ consists of $t_0'$. 
So every edge of the torso at $Q$ is either between $A$ and $B$, or an edge of the torso at $L$ with both ends in $A$.
Since there are at least $\lvert [A,B] \rvert+1$ edges of the torso at $L$ incident with $B$, the number of edges of the torso at $Q$ is strictly smaller than the number of edges of the torso at $L$.

So we may assume that $L$ is not a $\theta$-cell.
Hence $L=t_0$ and $T_0=T$.
So there exists a pseudo-$r$-cell $f(Q)$ in $(T,\X)$ such that $Q \subseteq f(Q)'$ by Claim 1.
Since $Q$ contains $t_0'$, $f(Q)$ contains $L$. 
Since $L \in V(f(Q))$ and there are at least $\lvert [A,B] \rvert+1$ edges of the torso at $L$ incident with $B$, there are at least $\lvert [A,B] \rvert+1$ edges of the torso at $f(Q)$ incident with $B$.
Since $Q$ contains $t_0'$ and $Q \subseteq T'$, every edge of the torso at $Q$ either has both ends in $A$ or is between $A$ and $B$.
Note that every edge in the former case is an edge of the torso at $f(Q)$ not incident with $B$.
Hence the number of edges of the torso at $f(Q)$ is strictly more than the number of edges of the torso at $Q$.
Since $Q$ is an $r$-cell and $f(Q)$ is a pseudo-$r$-cell, $f(Q)$ is an $r$-cell.
So $f(Q)$ is an $r$-cell containing $L$ with $k^* \leq r  \leq \theta$.
By the maximality of $k^*$, $k^* \geq r$.
So $k^*=r$ and hence $f(Q)=L^*$.
$\Box$

For every $r \in [\theta]$ with $r \geq k^*$, let $g_r$ be the function such that for every $r$-cell $Q$ of $(T^*,\X^*)$,
	\begin{itemize}
		\item if $Q$ contains $t_0'$ or $t_0''$, then $g_r(Q)=L^*$,
		\item otherwise, $g_r(Q)=f(Q)$.
	\end{itemize}
By Claims 2 and 3, for each $r$ with $\theta \geq r \geq k^*$, $g_r$ maps each $r$-cell $Q$ of $(T^*,\X^*)$ to an $r$-cell $g_r(Q)$ of $(T,\X)$ such that the number of edges of the torso at $Q$ is at most the number of edges of the torso at $g_r(Q)$; furthermore, if the number of edges of the torso at $Q$ and $g_r(Q)$ are the same, then there exists no $r$-cell $W$ other than $Q$ with $g_r(W)=g_r(Q)$.
Hence the $r$-signature of $(T^*,\X^*)$ is (lexicographically) at most the $r$-signature of $(T,\X)$ for each $r$ with $\theta \geq r \geq k^*$.
Since $L^*$ is a $k^*$-cell such that either there exists no $k^*$-cell in $(T^*,\X^*)$ mapped to $L^*$ by $g_{k^*}$, or all $k^*$-cells $Q$ of $(T^*,\X^*)$ with $g_{k^*}(Q)=L^*$ satisfy that the torso at $Q$ has less edges than the torso at $L^*$ by Claims 2 and 3, we know that the $\theta$-signature of $(T^*,\X^*)$ is lexicographically strictly smaller than the $\theta$-signature of $(T,\X)$, contradicting the minimality of $(T,\X)$.
This proves the theorem.
\end{pf}

\bigskip

Let $\theta$ be a positive integer.
An {\it edge-tangle} $\E$ of {\it order} $\theta$ in a graph $G$ is a set of edge-cuts of $G$ of order less than $\theta$ such that the following hold.
\begin{enumerate}
	\item[(E1)] For every edge-cut $[A,B]$ of $G$ of order less than $\theta$, either $[A,B] \in \E$ or $[B,A] \in \E$;
	\item[(E2)] If $[A_1,B_1],[A_2,B_2],[A_3,B_3] \in \E$, then $B_1 \cap B_2 \cap B_3 \neq \emptyset$.
	\item[(E3)] If $[A,B] \in \E$, then $G$ has at least $\theta$ edges incident with vertices in $B$.
\end{enumerate}
Note that if an edge-tangle $\E$ of order $\theta \geq 1$ in $G$ exists, then $[\emptyset, V(G)] \in \E$ by (E1) and (E2), so $\lvert E(G) \rvert \geq \theta$ by (E3).
Furthermore, for every $[A,B] \in \E$, there exists an edge of $G$ whose every end is in $B$ by (E3).

The following lemma (Lemma \ref{fat_cell_edge-tangle}) shows that every $\theta$-cell in a $\theta$-smooth tree-cut decomposition with sufficiently many edges in its torso defines an edge-tangle of order $\theta$.

Let $(T,\X)$ be a tree-cut decomposition of a graph $G$.
Let $\theta$ be a positive integer.
We say that a $\theta$-cell $L$ is {\it $\theta$-fat} if the torso at $L$ contains at least $3\theta-2$ edges.
Lemma \ref{fat_cell_edge-tangle} shows that every $\theta$-fat $\theta$-cell in a $\theta$-smooth tree-cut decomposition defines an edge-tangle of order $\theta$.
We call the edge-tangle $\E$ mentioned in Lemma \ref{fat_cell_edge-tangle} the {\it edge-tangle defined by $L$}.

\begin{lemma} \label{fat_cell_edge-tangle}
Let $G$ be a graph.
Let $\theta$ be a positive integer.
Let $(T,\X)$ be a $\theta$-smooth tree-cut decomposition in $G$.
Let $L$ be a $\theta$-fat $\theta$-cell.
Let $\E$ be the collection of edge-cuts $[A,B]$ of order less than $\theta$ such that $A$ is incident with at most $\lvert [A,B] \rvert$ edges of the torso at $L$.
Then $\E$ is an edge-tangle in $G$ of order $\theta$.
\end{lemma}

\begin{pf}
Let $[A,B]$ be an edge-cut of $G$ of order less than $\theta$.
If $[A,B] \not \in \E$ and $[B,A] \not \in \E$, then $A$ is incident with at least $\lvert [A,B] \rvert+1$ edges of the torso at $L$, and $B$ is incident with at least $\lvert [A,B] \rvert+1$ edges of the torso at $L$.
But $\lvert [A,B] \rvert < \theta$.
It contradicts that $(T,\X)$ is $\theta$-smooth.
Hence either $[A,B] \in \E$ or $[B,A] \in \E$.
Therefore, $\E$ satisfies (E1).

Suppose that there exist $[A_i,B_i] \in \E$ for $i \in [3]$ such that $B_1 \cap B_2 \cap B_3 = \emptyset$.
For each $i$, since $[A_i,B_i] \in \E$, $A_i$ is incident with at most $\lvert [A,B] \rvert \leq \theta-1$ edges of the torso at $L$.
Hence $A_1 \cup A_2 \cup A_3$ is incident with at most $3\theta-3$ edges of the torso at $L$.
Since the torso at $L$ has at least $3\theta-2$ edges, some edge of the torso at $L$ has all ends in $B_1 \cap B_2 \cap B_3$.
So $B_1 \cap B_2 \cap B_3 \neq \emptyset$.
Therefore, $\E$ satisfies (E2).

Suppose that there exists $[C,D] \in \E$ such that $G$ has less than $\theta$ edges incident with $D$.
Then $D$ is incident with at most $\theta-1$ edges of the torso at $L$.
Since $[C,D] \in \E$, $C$ is incident with at most $\theta-1$ edges of the torso at $L$.
Hence the torso at $L$ has at most $2\theta-2<3\theta-2$ edges, a contradiction.
Therefore, $\E$ satisfies (E3) and hence is an edge-tangle of order $\theta$
\end{pf}

\bigskip

The following simple lemma was proved in \cite{l_ep_imm}.

\begin{lemma}[{\cite[Lemma 2.3]{l_ep_imm}}] \label{easy_tangle}
Let $\theta$ be a positive integer.
Let $G$ be a graph and $\E$ an edge-tangle of order $\theta$ in $G$.
If $[A,B],[C,D] \in \E$, then the following hold.
	\begin{enumerate}
		\item If the order of $[A \cup C, B \cap D]$ is less than $\theta$, then $[A \cup C, B \cap D] \in \E$.
		\item If $A' \subseteq A$ and $[A',V(G)-A']$ is an edge-cut of $G$ of order less than $\theta$, then $[A',V(G)-A'] \in \E$.
	\end{enumerate}
\end{lemma}

The following lemma shows that given an edge-tangle $\E$ of a graph $G$, one can check, for each edge-cut $[A,B]$ of $G$ of small order given by an edge $e$ of a tree-cut decomposition $(T,\X)$, whether $[A,B] \in \E$, by simply seeing which component of $T-e$ contains a specific cell.

\begin{lemma} \label{weak_tangle_cell}
Let $G$ be a graph.
Let $\theta$ be a positive integer.
Let $\E$ be an edge-tangle in $G$ of order $\theta$.
Let $(T,\X)$ be a tree-cut decomposition of $G$. 
Then there exists a unique $\theta$-cell $C_\E$ such that for every edge $e$ of $T$ with at most one end in $V(C_\E)$ and with $\lvert \adh_{(T,\X)}(e) \rvert < \theta$, $[A_{e,C_\E},B_{e,C_\E}] \in \E$.
In particular, for every edge $e$ of $T$ with exactly one end in $V(C_\E)$, $[A_{e,C_\E},B_{e,C_\E}] \in \E$.
\end{lemma}

\begin{pf}
Note that for every edge $e=t_1t_2$ of $T$ with $\lvert \adh_{(T,\X)}(e) \rvert <\theta$, $[X_{V(T_{t_1,t_2})},X_{V(T_{t_2,t_1})}]$ is an edge-cut of $G$ of order less than $\theta$, where $T_{t_1,t_2},T_{t_2,t_1}$ are the two components of $T-e$ such that $T_{t_1,t_2}$ contains $t_2$, so either $[X_{V(T_{t_1,t_2})},X_{V(T_{t_2,t_1})}] \in \E$ or $[X_{V(T_{t_2,t_1})},X_{V(T_{t_1,t_2})}] \in \E$ by (E1).
Let $T^*$ be the tree obtained from $T$ by contracting each pseudo-$\theta$-cell into a node.
So for each $e \in E(T^*)$, $\lvert \adh_{(T,\X)}(e) \rvert < \theta$.
For every edge $t_1t_2$ of $T^*$, we define a direction on $t_1t_2$ such that $t_2$ is the head if and only if $[X_{V(T_{t_2,t_1})},X_{V(T_{t_1,t_2})}] \in \E$.
Note that the direction on the edges is well-defined by (E1) and (E2).
Since the sum of the out-degree of the nodes of $T^*$ equals $\lvert E(T^*) \rvert$, some node of $T^*$ has out-degree at most $\lvert E(T^*) \rvert/\lvert V(T^*) \rvert <1$.
So $T^*$ contains a node $t_\E$ with out-degree 0 in the orientation.

Suppose there exists an edge $e$ of $T^*$ such that $[A_{e,t_\E},B_{e,t_\E}] \not \in \E$.
Since $e$ is an edge of $T^*$, $\lvert \adh_{(T,\X)}(e) \rvert < \theta$.
Since $[A_{e,t_\E},B_{e,t_\E}] \not \in \E$ and $\lvert \adh_{(T,\X)}(e) \rvert < \theta$, by (E1), $[B_{e,t_\E},A_{e,t_\E}] \in \E$.
So $t_\E$ and the head of $e$ belong to different components of $T^*-e$.
Hence there exists an edge $e'$ of $T^*$ such that $e'$ is incident with $t_\E$, and $t_\E$ and the head of $e$ belong to different components of $T^*-e'$.
Since $t_\E$ has out-degree 0, $[A_{e',t_\E},B_{e',t_\E}] \in \E$.
Since $[B_{e,t_\E},A_{e,t_\E}] \in \E$ and $[A_{e',t_\E},B_{e',t_\E}] \in \E$, by (E2), $A_{e,t_\E} \cap B_{e',t_\E} \neq \emptyset$, a contradiction.

Hence for every edge $e$ of $T^*$, $[A_{e,t_\E},B_{e,t_\E}] \in \E$.
This implies that $t_\E$ is the unique node of $T^*$ with out-degree 0 in the orientation.

Note that there exists a pseudo-$\theta$-cell $T_\E$ in $(T,\X)$ contracted into $t_\E$.
Since every edge of $T$ with $\lvert \adh_{(T,\X)}(e) \rvert < \theta$ is an edge of $T^*$, for every edge $e$ of $T$ with at most one end in $T_\E$ and with $\lvert \adh_{(T,\X)}(e) \rvert < \theta$, $[A_{e,t_\E},B_{e,t_\E}] \in \E$.
Since $t_\E$ is the unique node of $T^*$ with out-degree 0, $T_\E$ is the unique pseudo-$\theta$-cell in $(T,\X)$ with this property.
Note that every edge $e$ of $T$ with exactly one end in $V(T_\E)$ has $\lvert \adh_{(T,\X)}(e) \rvert < \theta$ since $T_\E$ is a pseudo-$\theta$-cell.

To prove this lemma, it suffices to prove that $T_\E$ is a $\theta$-cell in $(T,\X)$.
That is, it suffices to prove that the torso at $T_\E$ in $(T,\X)$ contains at least $\theta$ edges.

Suppose that there are at most $\theta-1$ edges in the torso at $T_\E$ in $(T,\X)$.
Let $A_1 = \bigcup_e A_{e,T_\E}$, where the union is over all edges of $T$ with exactly one end in $T_\E$.
Let $B_1 = V(G)-A_1$.
So $B_1=\bigcap_e B_{e,T_\E}$, where the intersection is over all edges of $T$ with exactly one end in $T_\E$.
For every subset $S$ of the set of edges of $T$ with exactly one end in $T_\E$, let $A_S = \bigcup_{e \in S}A_{e,T_\E}$ and $B_S = V(G)-A_S$.
Note that for every such set $S$, $\lvert [A_S,B_S] \rvert$ is at most the number of edges of the torso at $T_\E$ which is less than $\theta$.
Since $[A_{e,T_\E}, B_{e,T_\E}] \in \E$ for each edge $e$ of $T$ with exactly one end in $T_\E$, by Lemma \ref{easy_tangle} and induction on $\lvert S \rvert$, we know $[A_1,B_1] \in \E$.

Note that $B_1=X_{V(T_\E)}$, and the number of edges of $G$ incident with $X_{V(T_\E)}$ is at most the number of edges of the torso at $T_\E$.
Hence $B_1$ is incident with at most $\theta-1$ edges of $G$.
So $[A_1,B_1] \not \in \E$ by (E3), a contradiction.
This shows that $T_\E$ is a $\theta$-cell.
\end{pf}

\bigskip

Let $\E_1,\E_2$ be distinct edge-tangles in a graph $G$. 
Then an {\it $(\E_1,\E_2)$-separator} is an edge-cut $[A,B]$ in $G$ such that $[A,B] \in \E_1-\E_2$ and $[B,A] \in \E_2-\E_1$.
A {\it minimum $(\E_1,\E_2)$-separator} is an $(\E_1,\E_2)$-separator with minimum order.

Let $G$ be a graph.
Let $\C$ be a collection of edge-tangles in $G$ and let $\E$ be an edge-tangle in $G$.
A {\it $(\C,\E)$-separator} is an edge-cut $[A,B]$ of $G$ such that $[A,B] \in \E'-\E$ and $[B,A] \in \E-\E'$ for every $\E' \in \C$.
A {\it minimum $(\C,\E)$-separator} is a $(\C,\E)$-separator with minimum order.

\begin{lemma} \label{simple_side}
Let $G$ be a graph.
Let $\C$ be a collection of edge-tangles in $G$.
Let $\E$ be an edge-tangle in $G$.
Let $[A,B]$ be a minimum $(\C,\E)$-separator.
Then for every $\E' \in \C$, there exists a minimum $(\E',\E)$-separator $[C,D]$ such that $A \subseteq C$ and $B \supseteq D$.
\end{lemma}

\begin{pf}
Let $\E'$ be a member of $\C$.
Let $[C,D]$ be a minimum $(\E',\E)$-separator such that $C$ is maximal.
Since $[A,B]$ is a $(\C,\E)$-separator, $[A,B]$ is an $(\E',\E)$-separator.

Suppose $\lvert [A \cup C, B \cap D] \rvert > \lvert [C,D] \rvert$.
By the submodularity, $\lvert [A \cap C,B \cup D] \rvert < \lvert [A,B] \rvert$.
Since $[A,B]$ is a $(\C,\E)$-separator, $[A \cap C, B \cup D]$ is a $(\C,\E)$-separator by Lemma \ref{easy_tangle}.
But the order of $[A \cap C, B \cup D]$ is smaller than $\lvert [A,B] \rvert$, a contradiction.

So $\lvert [A \cup C, B \cap D] \rvert \leq \lvert [C,D] \rvert$.
By Lemma \ref{easy_tangle}, $[A \cup C, B \cap D]$ is an $(\E',\E)$-separator of order at most $\lvert [C,D] \rvert$.
By the choice of $[C,D]$, $A \cup C \subseteq C$.
So $A \subseteq C$.
This implies that $B \supseteq D$.
\end{pf}

\bigskip

The following lemma builds a relationship between edge-tangles and fat cells in a smooth tree-cut decomposition and shows that information about edge-cuts separating different edge-tangles can be told by the tree edges.

\begin{lemma} \label{edge-tangle_tree}
Let $G$ be a connected graph.
Let $\theta$ be a positive integer.
Let $(T,\X)$ be a $\theta$-smooth tree-cut decomposition of $G$.
Then there exists a function $\iota$ that maps each edge-tangle $\E$ of order $\theta$ to a $\theta$-cell in $(T,\X)$ such that the following hold. 
	\begin{enumerate}
		\item For every edge-tangle $\E$ of order $\theta$ and each edge $e$ of $T$ with exactly one end in $V(\iota(\E))$, $[A_{e,\iota(\E)},B_{e,\iota(\E)}] \in \E$.
		\item If $L$ is a $\theta$-fat $\theta$-cell and $\E$ is an edge-tangle of order $\theta$ defined by $L$, then $\iota(\E)=L$.
		\item Let $\C$ be a collection of edge-tangles in $G$ of order $\theta$ such that every member of $\C$ is defined by some $\theta$-fat $\theta$-cell in $(T,\X)$.
		Let $\E$ be an edge-tangle in $G$ of order $\theta$ defined by some $\theta$-fat $\theta$-cell in $(T,\X)$.
		If there exists an edge $e$ of $T$ with exactly one end in $V(\iota(\E))$ such that the component of $T-e$ containing $\iota(\E)$ is disjoint from $V(\iota(\E'))$ for every $\E' \in \C$, then $[B_{e,\iota(\E)},A_{e,\iota(\E)}]$ is a $(\C,\E)$-separator, and for every minimum $(\C,\E)$-separator $[A,B]$, 
			\begin{enumerate}
				\item $B_{e,\iota(\E)} \subseteq A$, and
				\item for every $\E' \in \C$, if $e_{\E'}$ is the edge of $T$ with exactly one end in $V(\iota(\E'))$ such that $\iota(\E)$ and $\iota(\E')$ belong to different components of $T-e_{\E'}$, then $B_{e_{\E'},\iota(\E')} \subseteq B$.
			\end{enumerate}
	\end{enumerate}
\end{lemma}

\begin{pf}
Define $\iota$ to be the function whose domain is the set of all edge-tangles in $G$ of order $\theta$ such that for each edge-tangle $\E$, $\iota(\E)$ is the $\theta$-cell $C_\E$ mentioned in Lemma \ref{weak_tangle_cell}.
Then Statement 1 immediately follows from Lemma \ref{weak_tangle_cell}.
And Statement 2 follows from Lemma \ref{fat_cell_edge-tangle} and the uniqueness part of Lemma \ref{weak_tangle_cell}.

Now we prove Statement 3.
Let $\C,\E$ and $e$ be the collection, edge-tangle, and edge of $T$ mentioned in Statement 3, respectively.
By Statement 1, $[A_{e,\iota(\E)},B_{e,\iota(\E)}] \in \E$.
Since $e$ has exactly one end in $V(\iota(\E))$, $e$ has at most one end in $V(\iota(\E'))$ for every $\E' \in \C$.
For every $\E' \in \C$, since $\iota(\E')$ and $\iota(\E)$ are contained in different components of $T-e$, we have $[B_{e,\iota(\E)},A_{e,\iota(\E)}] = [A_{e,\iota(\E')},B_{e,\iota(\E')}] \in \E'$ by Lemma \ref{weak_tangle_cell}.
So $[B_{e,\iota(\E)},A_{e,\iota(\E)}]$ is a $(\C,\E)$-separator.

Let $[A,B]$ be a minimum $(\C,\E)$-separator. 
We first prove Statement 3(a).

If $\lvert [A \cup B_{e,\iota(\E)}, B \cap A_{e,\iota(\E)}] \rvert < \lvert [A,B] \rvert$, then by Lemma \ref{easy_tangle}, $[A \cup B_{e,\iota(\E)}, B \cap A_{e,\iota(\E)}]$ is a $(\C,\E)$-separator with order less than $\lvert [A,B] \rvert$, contradicting the minimality of $[A,B]$.
So $\lvert [A \cup B_{e,\iota(\E)}, B \cap A_{e,\iota(\E)}] \rvert \geq \lvert [A,B] \rvert$.
By the submodularity, $\lvert [A \cap B_{e,\iota(\E)}, B \cup A_{e,\iota(\E)}] \rvert \leq \lvert [A_{e,\iota(\E)},B_{e,\iota(\E)}] \rvert < \theta$.
So $[B \cup A_{e,\iota(\E)},A \cap B_{e,\iota(\E)}] \in \E$ and $[A \cap B_{e,\iota(\E)}, B \cup A_{e,\iota(\E)}]$ is a $(\C,\E)$-separator by Lemma \ref{easy_tangle}.

Let $G_1$ be the torso at $\iota(\E)$.
Since $\E$ is defined by $\iota(\E)$, $B \cup A_{e,\iota(\E)}$ is incident with at most $\lvert [B \cup A_{e,\iota(\E)},A \cap B_{e,\iota(\E)}] \rvert \leq \lvert [A_{e,\iota(\E)},B_{e,\iota(\E)}] \rvert$ edges of $G_1$.
Let $w_1$ be the vertex in $G_1$ obtained by identifying $A_{e,\iota(\E)}$.
Since every edge of $G_1$ incident with $w_1$ is an edge of $G$ incident with $B \cup A_{e,\iota(\E)}$, and there are $\lvert [A_{e,\iota(\E)},B_{e,\iota(\E)}] \rvert$ such edges, we know there are at least $\lvert [A_{e,\iota(\E)},B_{e,\iota(\E)}] \rvert$ edges of $G_1$ incident with $B \cup A_{e,\iota(\E)}$.
So $\lvert [B \cup A_{e,\iota(\E)},A \cap B_{e,\iota(\E)}] \rvert = \lvert [A_{e,\iota(\E)},B_{e,\iota(\E)}] \rvert$, and the set of edges of $G_1$ incident with $B \cup A_{e,\iota(\E)}$ equals the set of edges of $G_1$ between $A_{e,\iota(\E)}$ and $B_{e,\iota(\E)}$.
Since every edge of $G_1$ between $B \cup A_{e,\iota(\E)}$ and $A \cap B_{e,\iota(\E)}$ is incident with $B \cup A_{e,\iota(\E)}$, the set of edges of $G_1$ between $B \cup A_{e,\iota(\E)}$ and $A \cap B_{e,\iota(\E)}$ equals the set of edges of $G_1$ between $A_{e,\iota(\E)}$ and $B_{e,\iota(\E)}$.
In addition, since $\lvert [B \cup A_{e,\iota(\E)},A \cap B_{e,\iota(\E)}] \rvert = \lvert [A_{e,\iota(\E)},B_{e,\iota(\E)}] \rvert$ and every edge of $G$ between $A_{e,\iota(\E)}$ and $B_{e,\iota(\E)}$ is in $G_1$, the set of edges of $G$ between $B \cup A_{e,\iota(\E)}$ and $A \cap B_{e,\iota(\E)}$ equals the set of edges of $G$ between $A_{e,\iota(\E)}$ and $B_{e,\iota(\E)}$.

Suppose that there exists an edge $f$ of $G$ between $B_{e,\iota(\E)} \cap B$ and $A_{e,\iota(\E)} \cup A$.
Let $u$ be the end of $f$ in $B_{e,\iota(\E)} \cap B$.
Let $v$ be the end of $f$ in $A_{e,\iota(\E)} \cup A$.

We first assume that $v \in A_{e,\iota(\E)}$.
Then $f$ is between $A_{e,\iota(\E)}$ and $B_{e,\iota(\E)}$.
Recall that the set of edges of $G$ between $B \cup A_{e,\iota(\E)}$ and $A \cap B_{e,\iota(\E)}$ equals the set of edges of $G$ between $A_{e,\iota(\E)}$ and $B_{e,\iota(\E)}$.
So $f$ is between $B \cup A_{e,\iota(\E)}$ and $A \cap B_{e,\iota(\E)}$.
Since $u \in B$, $u \not \in A \cap B_{e,\iota(\E)}$.
Since $v \in A_{e,\iota(\E)}$, $v \not \in A \cap B_{e,\iota(\E)}$.
So no end of $f$ is in $A \cap B_{e,\iota(\E)}$, a contradiction.

Hence $v \in A-A_{e,\iota(\E)} = A \cap B_{e,\iota(\E)}$.
So $f$ is between $u \in B \subseteq B \cup A_{e,\iota(\E)}$ and $v \in A \cap B_{e,\iota(\E)}$.
Hence $f$ is between $A_{e,\iota(\E)}$ and $B_{e,\iota(\E)}$.
But both $u,v$ are contained in $B_{e,\iota(\E)}$, a contradiction.

Therefore, there exists no edge between $B_{e,\iota(\E)} \cap B$ and $A_{e,\iota(\E)} \cup A$.
Since $G$ is connected, either $B_{e,\iota(\E)} \cap B = \emptyset$ or $A_{e,\iota(\E)} \cup A = \emptyset$.
Since $[B,A] \in \E$, by (E2), $A \neq \emptyset$, so $A_{e,\iota(\E)} \cup A \neq \emptyset$.
Hence $B_{e,\iota(\E)} \cap B = \emptyset$.
That is, $B_{e,\iota(\E)} \subseteq A$.
This proves Statement 3(a).

Now we prove Statement 3(b).
Let $\E' \in \C$.
Since $[A,B]$ is a minimum $(\C,\E)$-separator, by Lemma \ref{simple_side}, there exists a minimum $(\E',\E)$-separator $[C,D]$ such that $B \supseteq D$.
Hence $[D,C]$ is a minimum $(\{\E\},\E')$-separator.
By applying Statement 3(a) with $\C=\{\E\}$ and $\E'=\E$, we know $B_{e_{\E'},\iota(\E')} \subseteq D$.
Hence $B_{e_{\E'},\iota(\E')} \subseteq D \subseteq B$.
This proves Statement 3(b).
\end{pf}

\section{Cross-free families} \label{sec:cross-free}

Let $\C$ and $\C'$ be collections of edge-tangles in a graph $G$.
A {\it $(\C,\C')$-segregator} is a set $\Se$ of edge-cuts of $G$ such that
	\begin{itemize}
		\item every member $[A,B]$ of $\Se$ is a minimum $(\C,\E')$-separator for some $\E' \in \C'$, and
		\item for every edge-tangle $\E'$ in $\C'$, there exists $[A,B] \in \Se$ such that either $[A,B]$ is a minimum $(\C,\E')$-separator, or $A' \subseteq A$ for some minimum $(\C,\E')$-separator $[A',B']$.
	\end{itemize}

A family $\D$ of edge-cuts of a graph is {\it cross-free} if $A \cap C = \emptyset$ for every pair of distinct edge-cuts $[A,B],[C,D]$ in $\D$.

\begin{lemma} \label{uncross_segregator}
Let $G$ be a graph.
Let $\C$ and $\C'$ be collections of edge-tangles in $G$. 
Let $\Se$ be a $(\C,\C')$-segregator. 
Then there exists a $(\C,\C')$-segregator $\Se^*$ such that $\bigcup_{[A,B] \in \Se}A = \bigcup_{[A,B] \in \Se^*}A$ and $\Se^*$ is a cross-free family.
\end{lemma}

\begin{pf}
Let $\Se^*$ be a $(\C,\C')$-segregator such that 
	\begin{itemize}
		\item[(i)] $\bigcup_{[A,B] \in \Se^*}A = \bigcup_{[A,B] \in \Se}A$, and
		\item[(ii)] subject to (i), $\sum_{[A,B] \in \Se^*}\lvert A \rvert$ is as small as possible. 
	\end{itemize}
Note that $\Se$ is a $(\C,\C')$-segregator, so $\Se^*$ exists.

There do not exist distinct members $[A_1,B_1]$ and $[A_2,B_2]$ of $\Se^*$ such that $A_1 \subseteq A_2$, for otherwise, $\Se^*-\{[A_1,B_1]\}$ is a $(\C,\C')$-segregator satisfying (i) but violating (ii). 

We shall prove that $\Se^*$ is a cross-free family.
Suppose to the contrary that there exist $[A,B] \in \Se^*$ and $[C,D] \in \Se^*$ such that $A \cap C \neq \emptyset$.
Note that $A-C \neq \emptyset \neq C-A$.
Since $\Se^*$ is a $(\C,\C')$-segregator, there exist $\E_1,\E_2 \in \C'$ such that $[A,B]$ is a minimum $(\C,\E_1)$-separator and $[C,D]$ is a minimum $(\C,\E_2)$-separator.

\noindent{\bf Claim 1:} $\lvert [A \cap C, B \cup D] \rvert < \min\{\lvert [A,B] \rvert, \lvert [C,D] \rvert\}$. 

\noindent{\bf Proof of Claim 1:}
If $\lvert [A \cup C, B \cap D] \rvert < \lvert [A,B] \rvert$, then by Lemma \ref{easy_tangle}, $[A \cup C, B \cap D]$ is a $(\C,\E_1)$-separator with order smaller than a minimum $(\C,\E_1)$-separator, a contradiction.
So $\lvert [A \cup C, B \cap D] \rvert \geq \lvert [A,B] \rvert$.
If $\lvert [A \cup C, B \cap D] \rvert = \lvert [A,B] \rvert$, then $[A \cup C, B \cap D]$ is a minimum $(\C,\E_1)$-separator by Lemma \ref{easy_tangle}, and $(\Se^*-\{[A,B],[C,D]\}) \cup \{[A \cup C, B \cap D]\}$ is a $(\C,\C')$-segregator satisfying (i) but violating (ii).
Hence $\lvert [A \cup C, B \cap D] \rvert > \lvert [A,B] \rvert$.
By the submodularity, $\lvert [A \cap C, B \cup D] \rvert < \lvert [C,D] \rvert$.
Similarly, $\lvert [A \cap C, B \cup D] \rvert < \lvert [A,B] \rvert$.
$\Box$

Let $\C_1 = \{\E \in \C': [A,B]$ is a minimum $(\C,\E)$-separator$\}$.
Let $\C_2 = \{\E \in \C': [C,D]$ is a minimum $(\C,\E)$-separator$\}$.
Note that $\E_1 \in \C_1$ and $\E_2 \in \C_2$.

\noindent{\bf Claim 2:} For every $\E \in \C_1 \cup \C_2$, $[A \cap C, B \cup D] \in \E$.

\noindent{\bf Proof of Claim 2:}
If there exists $\E \in \C_1$ such that $[A \cap C, B \cup D] \not \in \E$, then by Lemma \ref{easy_tangle} and Claim 1, $[A \cap C, B \cup D]$ is a $(\C,\E)$-separator with order smaller than $\lvert [A,B] \rvert$, contradicting that $[A,B]$ is a minimum $(\C,\E)$-separator.
So for every $\E \in \C_1$, $[A \cap C, B \cup D] \in \E$.
Similarly, for every $\E \in \C_2$, $[A \cap C, B \cup D] \in \E$.
$\Box$

In particular, $[A \cap C, B \cup D] \in \E_1$.
So if $\lvert [A \cap D, B \cup C] \rvert < \lvert [A,B] \rvert$, then $[A \cap D, B \cup C] \not \in \E_1$ (for otherwise, $[B,A],[A \cap C, B \cup D],[A \cap D, B \cup C]$ are members of $\E_1$ with $A \cap (B \cup D) \cap (B \cup C)=\emptyset$, contradicting (E2)), and it is a $(\C,\E_1)$-separator with order smaller than $\lvert [A,B] \rvert$ by Lemma \ref{easy_tangle}, a contradiction.
Hence $\lvert [A \cap D, B \cup C] \rvert \geq \lvert [A,B] \rvert$.
By the submodularity, $\lvert [B \cap C, A \cup D] \rvert \leq \lvert [C,D] \rvert$.

\noindent{\bf Claim 3:} $[B \cap C, A \cup D]$ is a minimum $(\C,\E)$-separator for every $\E \in \C_2$.

\noindent{\bf Proof of Claim 3:}
By Claim 2, for every $\E \in \C_2$, $[A \cap C, B \cup D] \in \E$, so $[B \cap C, A \cup D] \not \in \E$ by (E2).
So $[B \cap C, A \cup D]$ is a $(\C,\E)$-separator for every $\E \in \C_2$ by Lemma \ref{easy_tangle}.
Since $\lvert [B \cap C, A \cup D] \rvert \leq \lvert [C,D] \rvert$, we know $\lvert [B \cap C, A \cup D] \rvert = \lvert [C,D] \rvert$ and $[B \cap C, A \cup D]$ is a minimum $(\C,\E)$-separator for every $\E \in \C_2$.
$\Box$

Let $\Se' = (\Se^*-\{[C,D]\}) \cup \{[B \cap C, A \cup D]\}$.
Then $\Se'$ satisfies (i).
By (ii), $\Se'$ is not a $(\C,\C')$-segregator.
Since $[B \cap C, A \cup D]$ is a minimum $(\C,\E_2)$-separator by Claim 3, there exists $\E^* \in \C'$ such that 
	\begin{itemize}
		\item[(a)] either $[C,D]$ is a minimum $(\C,\E^*)$-separator, or $C' \subseteq C$ for some minimum $(\C,\E^*)$-separator $[C',D']$, 
		\item[(b)] $[B \cap C, A \cup D]$ is not a minimum $(\C,\E^*)$-separator, and 
		\item[(c)] for every minimum $(\C,\E^*)$-separator $[C',D']$, $C' \not \subseteq B \cap C$ and $C' \not \subseteq A$.
	\end{itemize}

By (b) and Claim 3, $\E^* \not \in \C_2$.
So $[C,D]$ is not a minimum $(\C,\E^*)$-separator.
By (a), there exists a minimum $(\C,\E^*)$-separator $[C',D']$ such that $C' \subseteq C$.

If $\lvert [A \cup C', B \cap D'] \rvert < \lvert [A,B] \rvert$, then $[A \cup C', B \cap D']$ is a $(\C,\E_1)$-separator with order smaller than $\lvert [A,B] \rvert$ by Lemma \ref{easy_tangle}, contradicting that $[A,B]$ is a minimum $(\C,\E_1)$-separator.
So $\lvert [A \cup C', B \cap D'] \rvert \geq \lvert [A,B] \rvert$.
By the submodularity, $\lvert [A \cap C', B \cup D'] \rvert \leq \lvert [C',D'] \rvert$.

Suppose $[A \cap C', B \cup D'] \not \in \E^*$.
Since $\lvert [A \cap C', B \cup D'] \rvert \leq \lvert [C',D'] \rvert$, $[A \cap C', B \cup D']$ is a minimum $(\C,\E^*)$-separator by Lemma \ref{easy_tangle}.
By (c), $A \cap C' \not \subseteq A$, a contradiction.

Hence $[A \cap C', B \cup D'] \in \E^*$.
Suppose $\lvert [A \cup D', B \cap C'] \rvert \leq \lvert [C', D'] \rvert$.
Since $[D',C'] \in \E^*$ and $[A \cap C', B \cup D'] \in \E^*$, $[A \cup D', B \cap C'] \in \E^*$ by Lemma \ref{easy_tangle}.
By Lemma \ref{easy_tangle}, $[B \cap C', A \cup D']$ is a $(\C,\E^*)$-separator.
Since $\lvert [A \cup D', B \cap C'] \rvert \leq \lvert [C', D'] \rvert$, $[B \cap C', A \cup D']$ is a minimum $(\C,\E^*)$-separator.
By (c), $B \cap C' \not \subseteq B \cap C$.
However, $C' \subseteq C$, so $B \cap C' \subseteq B \cap C$, a contradiction.

Hence $\lvert [A \cup D', B \cap C'] \rvert > \lvert [C',D'] \rvert$.
By the submodularity, $\lvert [A \cap D', B \cup C'] \rvert < \lvert [A,B] \rvert$.
Since $[A,B] \in \E$ for every $\E \in \C$, $[A \cap D', B \cup C'] \in \E$ for every $\E \in \C$ by Lemma \ref{easy_tangle}.
Since $[A,B]$ is a minimum $(\C,\E_1)$-separator and $\lvert [A \cap D', B \cup C'] \rvert < \lvert [A,B] \rvert$, $[A \cap D', B \cup C'] \in \E_1$.
By Claim 2, $[B,A], [A \cap C, B \cup D], [A \cap D', B \cup C']$ are members of $\E_1$ such that $A \cap (B \cup D) \cap (B \cup C') = A \cap D \cap C' \subseteq A \cap D \cap C = \emptyset$, contradicting (E2).
This proves the lemma.
\end{pf}

\bigskip

Let $\theta$ be a positive integer.
Let $G$ be a graph.
Let $\D$ be a cross-free family in $G$.
Let $\C$ be a collection of edge-tangles in $G$ of order $\theta$ such that $\D \subseteq \E$ for every $\E \in \C$.
A {\it $(\D,\C)$-guard} is a cross-free family $\Se$ such that for every edge-cut $[C,D] \in \Se$, $[C,D]$ is a minimum $(\C,\E')$-separator for some edge-tangle $\E'$ of order $\theta$ with $\D \not \subseteq \E'$.

Intuitively, a $(\D,\C)$-guard is a collection of edge-cuts distinguishing all edge-tangles in $\C$ from some edge-tangles inconsistent with $\D$.
The following lemma shows that given a $(\D,\C)$-guard $\Se$, one can update $\D$ so that $\D$ ``simulates'' the edge-cuts in $\Se$ without losing some nice properties of $\D$.

\begin{lemma} \label{insert_guard}
Let $G$ be a graph.
Let $\theta,\xi,d$ be positive integers. 
Let $\D_0$ be a cross-free family such that there exists $\D_0' \subseteq \D_0$ with $\lvert \D_0' \rvert \leq \xi$ such that every edge-cut in $\D_0-\D_0'$ has order at most $d-1$, and every edge-cut in $\D_0'$ has order at most $d-1+\xi$. 
Let $\C$ be a collection of edge-tangles in $G$ of order $\theta$ such that $\D_0 \subseteq \E$ for every $\E \in \C$.
Let $\Se$ be a $(\D_0,\C)$-guard.
If $\theta \geq d+\xi$, then there exist cross-free families $\D$ and $\D'$ such that 
	\begin{enumerate}
		\item $\D' \subseteq \D$. 
		\item $\lvert \D' \rvert \leq \xi$, 
		\item every edge-cut in $\D-\D'$ has order at most $d-1$, and every edge-cut in $\D'$ has order at most $d-1+\xi$, 
		\item $\bigcup_{[A,B] \in \D}A \supseteq \bigcup_{[A,B] \in \D_0}A$, 
		\item $\D \subseteq \E$ for every $\E \in \C$, and 
		\item for every $[C,D] \in \Se$, there exists $[A,B] \in \D$ such that $C \subseteq A$. 
	\end{enumerate}
\end{lemma}

\begin{pf}
	For every cross-free family $\F$, define $\Se_{\F}$ to be the set $\{[C,D] \in \Se:$ there exists $[A,B] \in \F$ with $C \subseteq A\}$. 
	We say that a cross-free family $\F$ is {\it useful} if
	\begin{itemize}
		\item[(i)] there exists $\F' \subseteq \F$ with $\lvert \F' \rvert \leq \xi$ such that every edge-cut in $\F-\F'$ has order at most $d-1$, and every edge-cut $[A,B]$ in $\F'$ has order at most $d-1+\xi$, 
		\item[(ii)] $\F \subseteq \E$ for every $\E \in \C$, 
		\item[(iii)] $\bigcup_{[A,B] \in \F}A \supseteq \bigcup_{[A,B] \in \D_0}A$, and
		\item[(iv)] if $\E'$ is an edge-tangle of order $\theta$ in $G$ with $\D_0 \not \subseteq \E'$ such that some member of $\Se$ is a minimum $(\C,\E')$-separator, then $\F \not \subseteq \E'$.
	\end{itemize}
	We call $(\F,\F')$ a {\it witness} if $\F,\F'$ satisfy (i).
	Note that $\D_0$ is useful and $(\D_0,\D_0')$ is a witness.
	
	Let $\D$ be a useful cross-free family such that $\lvert \Se_{\D} \rvert$ is as large as possible. 
	Let $(\D,\D')$ be a witness.
	We shall prove that $\D$ and $\D'$ satisfy all conclusions of this lemma.
	Clearly $\D$ and $\D'$ satisfy Statements 1-5.
	
	Suppose that $\D$ does not satisfy Statement 6.
	So there exists $[C,D] \in \Se-\Se_\D$. 
	That is, for every $[A,B] \in \D$, $C \not \subseteq A$.
	
	Since $[C,D] \in \Se$ and $\Se$ is a $(\D_0,\C)$-guard, there exists an edge-tangle $\E'$ of order $\theta$ such that $\D_0 \not \subseteq \E'$ and $[C,D]$ is a minimum $(\C,\E')$-separator.
	By (iv), $\D \not \subseteq \E'$, so there exists $[A_{\E'},B_{\E'}] \in \D-\E'$.
	
	\noindent{\bf Claim 1:} $\lvert [C,D] \rvert \leq \lvert [A_{\E'},B_{\E'}] \rvert$, $[D,C] \in \E'$ and $[B_{\E'},A_{\E'}] \in \E'$.
	
	\noindent{\bf Proof of Claim 1:}
	Since $\theta>d-1+\xi \geq \lvert [A_{\E'},B_{\E'}] \rvert$, $[B_{\E'},A_{\E'}] \in \E'$ by (E1).
	Since $[A_{\E'},B_{\E'}] \in \D-\E'$, $[A_{\E'},B_{\E'}]$ is a $(\C,\E')$-separator by (ii).
	Since $[C,D]$ is a minimum $(\C,\E')$-separator, $\lvert [C,D] \rvert \leq \lvert [A_{\E'},B_{\E'}] \rvert$, so $[D,C] \in \E'$ by (E1).
	$\Box$
	
	Define $\Q = \{[A_{\E'} \cup C, B_{\E'} \cap D]\} \cup \{[A \cap D, B \cup C]: [A,B] \in \D-\{[A_{\E'},B_{\E'}]\}\}$.
	Let $\Q_0' = \{[A \cap D, B \cup C]: [A,B] \in \D'-\{[A_{\E'},B_{\E'}]\}\}$.
	If $[A_{\E'},B_{\E'}] \in \D'$, then define $\Q' = \Q_0' \cup \{[A_{\E'} \cup C, B_{\E'} \cap D]\}$; otherwise, define $\Q' = \Q_0'$.
	Since $\D$ is a cross-free family, $\Q$ is a cross-free family.
	
\noindent{\bf Claim 2:} $\Se_\Q \supseteq \Se_\D$.

\noindent{\bf Proof of Claim 2:}
Let $[X,Y] \in \Se_\D$.
So there exists $[X^*,Y^*] \in \D$ such that $X \subseteq X^*$. 
Since $[X,Y] \in \Se_\D$ and $[C,D] \not \in \Se_\D$, $[X,Y] \neq [C,D]$.
Since $\Se$ is cross-free and $\Se_\D \subseteq \Se$, $X \cap C=\emptyset$.
So $X \subseteq D$.
If $[X^*,Y^*]=[A_{\E'},B_{\E'}]$, then $X \subseteq X^* \subseteq A_{\E'} \cup C$, so $[X,Y] \in \Se_\Q$.
So we may assume $[X^*,Y^*] \neq [A_{\E'},B_{\E'}]$.
That is, $[X^*,Y^*] \in \D-\{[A_{\E'},B_{\E'}]\}$, and $[X^* \cap D, Y^* \cup C] \in \Q$.
Note that $X \subseteq X^* \cap D$.
So $[X,Y] \in \Se_\Q$.
$\Box$

	Since $\Se$ is a cross-free family and $[C,D] \in \Se \cap \Se_Q-\Se_\D$, $\Se_\Q \supset \Se_\D$ by Claim 2.
Hence by the choice of $\D$, $\Q$ is not useful.
	Clearly, $\Q$ satisfies (iii).
	Note that $\lvert \Q' \rvert \leq \lvert \D' \rvert \leq \xi$.
	
	\noindent{\bf Claim 3:} $\lvert [A_{\E'} \cup C, B_{\E'} \cap D] \rvert \leq \lvert [A_{\E'},B_{\E'}] \rvert$. 
	
	\noindent{\bf Proof of Claim 3:}
	Suppose $\lvert [A_{\E'} \cup C, B_{\E'} \cap D] \rvert > \lvert [A_{\E'},B_{\E'}] \rvert$.
	By the submodularity, $\lvert [A_{\E'} \cap C, B_{\E'} \cup D] \rvert < \lvert [C,D] \rvert$.
	Since $[C,D]$ is a $(\C,\E')$-separator, $[C,D] \in \E-\E'$ and $[D,C] \in \E'-\E$ for every $\E \in \C$.
	Since $[A_{\E'},B_{\E'}] \in \D-\E'$, $[A_{\E'},B_{\E'}] \in \E-\E'$ for every $\E \in \C$.
	Hence $[A_{\E'} \cap C, B_{\E'} \cup D] \in \E-\E'$ for every $\E \in \C$ by Lemma \ref{easy_tangle}.
	Since $\lvert [A_{\E'} \cap C, B_{\E'} \cup D] \rvert < \lvert [C,D] \rvert$, by (E1), $[B_{\E'} \cup D, A_{\E'} \cap C] \in \E'-\E$ for every $\E \in \C$.
	So $[A_{\E'} \cap C, B_{\E'} \cup D]$ is a $(\C,\E')$-separator with order smaller than $\lvert [C,D] \rvert$, a contradiction.
	Hence $\lvert [A_{\E'} \cup C, B_{\E'} \cap D] \rvert \leq \lvert [A_{\E'},B_{\E'}] \rvert$.
	$\Box$
	
	\noindent{\bf Claim 4:} For every $[A,B] \in \D-\{[A_{\E'},B_{\E'}]\}$, $\lvert [A \cap D, B \cup C] \rvert \leq \lvert [A,B] \rvert$.
	
	\noindent{\bf Proof of Claim 4:}
	Suppose to the contrary that there exists $[A,B] \in \D-\{[A_{\E'},B_{\E'}]\}$ such that $\lvert [A \cap D, B \cup C] \rvert > \lvert [A,B] \rvert$.
	By the submodularity, $\lvert [A \cup D, B \cap C] \rvert < \lvert [C,D] \rvert<\theta$.
	Since $[C,D] \in \E-\E'$ for every $\E \in \C$, by Lemma \ref{easy_tangle}, $[B \cap C, A \cup D] \in \E$ for every $\E \in \C$.
	Since $\D$ is a cross-free family, $A_{\E'} \cap A = \emptyset$.
	Since $[B_{\E'},A_{\E'}] \in \E'$ and $[D,C] \in \E'$ by Claim 1, $[B \cap C, A \cup D] \not \in \E'$ by (E2).
	Hence $[B \cap C, A \cup D]$ is a $(\C,\E')$-separator with order smaller than $\lvert [C,D] \rvert$, a contradiction.
	$\Box$
	
	By Claims 3 and 4, $(\Q,\Q')$ is a witness.

\noindent{\bf Claim 5:} If $\E''$ is an edge-tangle of order $\theta$ in $G$ with $\D_0 \not \subseteq \E''$ such that some member of $\Se$ is a minimum $(\C,\E'')$-separator, then $\Q \not \subseteq \E''$.

\noindent{\bf Proof of Claim 5:}
Suppose to the contrary that there exists an edge-tangle $\E''$ of order $\theta$ with $\D_0 \not \subseteq \E''$ and $\Q \subseteq \E''$ and there exists a minimum $(\C,\E'')$-separator $[X,Y] \in \Se$.
Since $[A_{\E'} \cup C, B_{\E'} \cap D] \in \Q \subseteq \E''$, by Lemma \ref{easy_tangle}, $[A_{\E'},B_{\E'}] \in \E''$ and $[C,D] \in \E''$.
Since $\D$ satisfies (iv), $\D \not \subseteq \E''$.
So there exists $[A,B] \in \D-\{[A_{\E'},B_{\E'}]\}$ such that $[A,B] \not \in \E''$.
Hence $[A \cap D, B \cup C] \in \Q \subseteq \E''$.
However, $[B,A], [C,D], [A \cap D, B \cup C]$ are members of $\E''$ with $A \cap D \cap (B \cup C)=\emptyset$, contradicting (E2).
$\Box$

Hence $\Q$ satisfies (iv) by Claim 5.
Recall that $\Q$ is not useful but satisfies (iii). 
So $\Q$ does not satisfy (ii). 
That is, there exists $\E^* \in \C$ such that $\Q \not \subseteq \E^*$. 
So there exists $[A^*,B^*] \in \Q-\E^*$.
	
Since $\E^* \in \C$, we know that $[A_{\E'},B_{\E'}] \in \D \subseteq \E^*$ and $[C,D] \in \E^*$.
So $[A_{\E'} \cup C, B_{\E'} \cap D] \in \E^*$ by Claim 3 and Lemma \ref{easy_tangle}.
Hence $[A^*,B^*] = [A' \cap D, B' \cup C]$ for some $[A',B'] \in \D-\{[A_{\E'},B_{\E'}]\}$.
	Since $[A^*,B^*] \not \in \E^*$, by Claim 4 and (E1), $[B' \cup C, A' \cap D]=[B^*,A^*] \in \E^*$.
	Since $\D$ satisfies (ii), $[A',B'] \in \E^*$.
	But $B' \cap (A' \cap D) \subseteq B' \cap A' = \emptyset$, contradicting (E2).
	This proves the lemma.
\end{pf}

\begin{lemma} \label{insert_one}
	Let $G$ be a graph.
	Let $\D$ be a cross-free family.
	Let $[A^*,B^*]$ be an edge-cut of $G$. 
	For every $[A,B] \in \D$, let $f([A,B])=[A \cap B^*,B \cup A^*]$.
	Let $\D' = \{[A^*,B^*]\} \cup \{f([A,B]): [A,B] \in \D\}$.
	Then there exists $\D'' \subseteq \D'-\{[A^*,B^*]\}$ with $\lvert \D'' \rvert \leq 2\lvert [A^*,B^*] \rvert$ such that for every $f([A,B]) \in \D'-(\D'' \cup \{[A^*,B^*]\})$, $\lvert f([A,B]) \rvert \leq \lvert [A,B] \rvert$, and for every $f([A,B]) \in \D''$, $\lvert f([A,B]) \rvert \leq \lvert [A,B] \rvert + \lvert [A^*,B^*] \rvert$.
\end{lemma}

\begin{pf}
	Define $\D''=\{f([A,B]):[A,B] \in \D$ and some edge between $A^*$ and $B^*$ is between $A \cap B^*$ and $A^*\}$.
	Note that for every $[A,B] \in \D$, if $\lvert f([A,B]) \rvert > \lvert [A,B] \rvert$, then some edge between $A^*$ and $B^*$ is between $A \cap B^*$ and $A^*$, so $f([A,B]) \in \D''$.
	Furthermore, for every $[A,B] \in \D$, $\lvert f([A,B]) \rvert \leq \lvert [A,B] \rvert + \lvert [A^*,B^*] \rvert$ by the subdomularity.
	In addition, since there are at most $\lvert [A^*,B^*] \rvert$ edges between $A^*$ and $B^*$, and for distinct $[A_1,B_1],[A_2,B_2] \in \D$, $A_1 \cap B^*$ and $A_2 \cap B^*$ are disjoint, we know $\lvert \D'' \rvert \leq 2\lvert [A^*,B^*] \rvert$.
	This proves the lemma.
\end{pf}

\section{Excluding immersions in 3-edge-cut free graphs} \label{sec:4-edge-conn} 

Now we focus on graphs with no edge-cut of order 3. 
The following theorem proved in \cite{l_ep_imm} states that in any graph with no edge-cut of order 3, every edge-tangle of large order ``controls a $K_k$-thorns'' for some large $k$.
(We omit the formal definition for ``an edge-tangle controlling a $K_{k}$-thorns'' mentioned in the following theorem, because it requires some sentences to be formally stated and we do not need this formal definition in this paper.)

\begin{lemma}[{\cite[Theorem 6.4]{l_ep_imm}}] \label{edge-tangle_4-edge-connected_control_edge-minor}
For any positive integers $k$ and $\theta$ with $\theta > k$, there exists a positive integer $w=w(k,\theta)$ such that if $G$ is a graph with no edge-cut of order 3, and $\E$ is an edge-tangle in $G$ of order at least $w$, then $\E_\theta$ controls a $K_k$-thorns, where $\E_\theta$ is the edge-tangle in $G$ of order $\theta$ such that $\E_\theta \subseteq \E$.
\end{lemma}

Let $G$ be a graph and $\E$ a collection of edge-cuts of $G$ of order less than a positive integer $\theta$, and let $X \subseteq E(G)$.
Define $\E-X$ to be the set of edge-cuts of $G-X$ of order less than $\theta-\lvert X \rvert$ such that $[A,B] \in \E-X$ if and only if $[A,B] \in \E$.
The following is proved in \cite{l_ep_imm}.

\begin{lemma}[{\cite[Lemma 2.6]{l_ep_imm}}] \label{edge-tangle_deleting_edges}
Let $G$ be a graph and $\theta$ be a positive integer.
If $\E$ is an edge-tangle in $G$ of order $\theta$ and $X$ is a subset of $E(G)$ with $\lvert X \rvert < \theta$, then $\E-X$ is an edge-tangle in $G-X$ of order $\theta-\lvert X \rvert$.
\end{lemma}

Recall that a graph is {\it exceptional} if it contains exactly one vertex of degree at least two, and this vertex is incident with a loop.
The following is the structure theorem for excluding a non-exceptional graph as an immersion with respect to an edge-tangle controlling a ``thorns'' proved in \cite{l_ep_imm}.

\begin{lemma}[{\cite[Theorem 4.6]{l_ep_imm}}] \label{excluding_immersion_thorns}
For any positive integers $d$ and $h$, there exist positive integers $\theta=\theta(d,h)$ and $\xi=\xi(d,h)$ such that the following holds.
If $H$ is a non-exceptional graph with degree sequence $(d_1,d_2,...,d_h)$, where $d_1=d$, and $G$ is a graph that does not contain an $H$-immersion, then for every edge-tangle $\E$ of order at least $\theta$ in $G$ controlling a $K_{3dh}$-thorns, there exist $C \subseteq E(G)$ with $\lvert C \rvert \leq \xi$, $U \subseteq V(G)$ with $\lvert U \rvert \leq h-1$ and a cross-free family ${\mathcal D} \subseteq \E-C$ such that for every vertex $v \in V(G)-U$, there exists $[A,B] \in {\mathcal D}$ of order at most $d_{\lvert U \rvert+1}-1$ with $v \in A$.
\end{lemma}

By simply combining Lemmas \ref{edge-tangle_4-edge-connected_control_edge-minor} and \ref{excluding_immersion_thorns}, we obtain the following result for excluding a non-exceptional graph as an immersion in graphs with no edge-cut of order 3 with respect to an edge-tangle.
This result will be used in proving Lemma \ref{global_no_iso_5}.

\begin{lemma} \label{basic_exclu_imm}
For any positive integers $d$ and $h$, there exist positive integers $\theta=\theta(d,h)$ and $\xi=\xi(d,h)$ such that the following holds.
If $H$ is a non-exceptional graph with degree sequence $(d_1,d_2,...,d_h)$, where $d_1=d$, and $G$ is a graph that does not contain an $H$-immersion and has no edge-cut of order 3, then for every edge-tangle $\E$ of order at least $\theta$ in $G$, there exist $C \subseteq E(G)$ with $\lvert C \rvert \leq \xi$, $U \subseteq V(G)$ with $\lvert U \rvert \leq h-1$ and a cross-free family ${\mathcal D} \subseteq \E-C$ such that for every vertex $v \in V(G)-U$, there exists $[A,B] \in {\mathcal D}$ of order at most $d_{\lvert U \rvert+1}-1$ with $v \in A$.
\end{lemma}

\begin{pf}
We define the following.
	\begin{itemize}
		\item Let $\theta_0 = \theta_{\ref{excluding_immersion_thorns}}(d,h)$, where $\theta_{\ref{excluding_immersion_thorns}}$ is the integer $\theta$ mentioned in Lemma \ref{excluding_immersion_thorns}.
		\item Let $\theta_1 = \theta_0+3dh$.
		\item Define $\theta = w_{\ref{edge-tangle_4-edge-connected_control_edge-minor}}(3dh,\theta_1)$, where $w_{\ref{edge-tangle_4-edge-connected_control_edge-minor}}$ is the integer $w$ mentioned in Lemma \ref{edge-tangle_4-edge-connected_control_edge-minor}.
		\item Define $\xi = \xi_{\ref{excluding_immersion_thorns}}(d,h)$, where $\xi_{\ref{excluding_immersion_thorns}}$ is the integer $\xi$ mentioned in Lemma \ref{excluding_immersion_thorns}.
	\end{itemize}

Let $H$ be a non-exceptional graph with degree sequence $(d_1,d_2,...,d_h)$ with $d_1=d$.
Let $G$ be a graph that does not contain an $H$-immersion and has no edge-cut of order 3.
Let $\E$ be an edge-tangle in $G$ of order at least $\theta$.
By Lemma \ref{edge-tangle_4-edge-connected_control_edge-minor}, $\E_{\theta_1}$ controls a $K_{3dh}$-thorns, where $\E_{\theta_1}$ is the edge-tangle in $G$ of order $\theta_1$ such that $\E_{\theta_1} \subseteq \E$.
By Lemma \ref{excluding_immersion_thorns}, there exist $C \subseteq E(G)$ with $\lvert C \rvert \leq \xi$, $U \subseteq V(G)$ with $\lvert U \rvert \leq h-1$ and a cross-free family $\D \subseteq \E_{\theta_1}-C \subseteq \E-C$ such that for every vertex $v \in V(G)-U$, there exists $[A,B] \in \D$ of order at most $d_{\lvert U \rvert+1}-1$ with $v \in A$.
This proves the lemma.
\end{pf}

\bigskip

A {\it rooted tree} $T$ rooted at a node $r$ is a directed graph whose underlying graph is a tree such that for every node $t$ of $T$, there exists a directed path in $T$ from $r$ to $t$.
A node $t_1$ of $T$ is an {\it ancestor} of a node $t_2$ of $T$ if there exists a directed path in $T$ from $t_1$ to $t_2$.
We say that $t_2$ is a {\it descendant} of $t_1$ if $t_1$ is an ancestor of $t_2$.
Note that every node is an ancestor and a descendant of itself.

Lemma \ref{global_no_iso_5} is the heart of our global decomposition theorem.
We sketch its proof here.
By Theorem \ref{theta-smooth} and Lemma \ref{fat_cell_edge-tangle}, for a sufficiently large $\theta$, we can obtain a $\theta$-smooth tree-cut decomposition $(T,\X)$ such that each $\theta$-fat $\theta$-cell defines an edge-tangle of order $\theta$.
For the simplicity of the proof sketch, we do not distinguish a $\theta$-fat $\theta$-cell and the edge-tangle defined by it.
For each edge-tangle $\E$ defined by a fat cell (call an ``important edge-tangle'' for simplicity), there exists a cross-free family $\D_\E$ that has a nice property as what the family $\D$ has in Lemma \ref{basic_exclu_imm}.
Then for each important edge-tangle $\E$, we define $\C_\E$ to be the maximal collection of important edge-tangles $\E'$ such that $\D_\E \subseteq \E'$ and those edge-tangles ``form a connected subtree'' in $T$, and define $\overline{\C_\E}$ to be the collection of the important edge-tangles not in $\C_\E$.
Note that we can define $\I$ to be a collection of important edge-tangles such that the ``connected subtrees in $T$ formed by $\C_{\E'}$'' for all $\E' \in \I$ form a partition of $T$.
We expect to define a new tree-cut decomposition $(T^*,\X^*)$ satisfying the conclusion of the lemma by first for each $\E \in \I$, replacing the ``connected subtree in $T$ formed by $\C_\E$'' together with the bag of the nodes contained in there by a new tree-cut decomposition $(T^\E,\X^\E)$ whose every node with sufficiently many edges in its torso``realizes'' $\D_\E$, and then by properly adding edges between the trees in $\{T^\E: \E \in \I\}$ to obtain $T^*$.
There are two potential concerns in this strategy. 
The first concern is that it is unclear how to add edges between those $T^\E$'s to obtain $T$.
This concern can be resolved by first using earlier lemmas to show that for each $\E \in \I$, there exists a $(\C_\E,\overline{\C_\E})$-segregator $\Se_\E$ that is also a $(\D_\E,\C_\E)$-guard, and then using the information of $\Se_\E$ to add edges to construct $T$.
The second concern is that it is unclear whether the resulting $(T^*,\X^*)$ satisfies the conclusion of the lemma, as those $\D_\E$'s can be very ``inconsistent'' with the edge-cuts given by $(T,\X)$.
To resolve this concern, we use earlier lemmas to show that for each important edge-tangle $\E$, we can modify $\D_\E$ to obtain another cross-free family $\D_{\E,2}$ that is more consistent with $(T,\X)$ than $\D_\E$ and still has nice properties so that if we use $\D_{\E,2}$ instead of $\D_\E$ to construction $(T^*,\X^*)$, then $(T^*,\X^*)$ satisfies the conclusion of this lemma.

\begin{lemma} \label{global_no_iso_5}
For any positive integers $d$ and $h$, there exist integers $\eta=\eta(d,h) \geq d$ and $\xi=\xi(d,h)$ such that the following holds.
Let $H$ be a graph with degree sequence $(d_1,d_2,...,d_h)$, where $d_1=d$, such that $H$ has no isolated vertex.
Let $G$ be a graph with no edge-cut of order 3 such that $G$ does not contain an $H$-immersion.
Define $H'=H$ if $H$ is non-exceptional; otherwise, define $H'$ to be a graph obtained from $H$ by subdividing one edge.
Let $h'=\lvert V(H') \rvert$ and $(d_1',d_2',...,d'_{h'})$ be the degree sequence of $H'$.
Then there exists a tree-cut decomposition $(T,\X=(X_t: t \in V(T)))$ of $G$ of adhesion at most $\eta$ such that for every $t \in V(T)$, there exist $Z_t \subseteq E(G)$ and $U_t \subseteq X_t$ with $\lvert Z_t \rvert \leq \xi$ and $\lvert U_t \rvert \leq h'-1$ such that if $G_t$ is the torso at $t$, then for every $v \in V(G_t-Z_t)-U_t$, the degree of $v$ in $G_t-Z_t$ is at most $d'_{\lvert U_t \rvert+1}-1$.
\end{lemma}

\begin{pf}
We define the following.
	\begin{itemize}
		\item Let $\theta_0=\max\{\theta_{\ref{basic_exclu_imm}}(d,h),\theta_{\ref{basic_exclu_imm}}(d,h+1)\}$ and $\xi_0=\max\{\xi_{\ref{basic_exclu_imm}}(d,h),\xi_{\ref{basic_exclu_imm}}(d,h+1)\}$, where $\theta_{\ref{basic_exclu_imm}}$ and $\xi_{\ref{basic_exclu_imm}}$ are the integers $\theta$ and $\xi$ mentioned in Lemma \ref{basic_exclu_imm}.
		\item Let $\theta = \theta_0+d+2\xi_0$.
		\item Define $\eta = 5\theta$.
		\item Define $\xi = 8\theta^2$.
	\end{itemize}

Let $H,G,H'$ be the graphs and let $h'$ be the integer as stated in the lemma.
Note that we may assume that $G$ is connected since if the lemma holds for the case when $G$ is connected, then when $G$ is disconnected, we can first obtain a desired tree-cut decomposition for each component of $G$, and then add edges between the underlying trees of those tree-cut decompositions for the components of $G$ to obtain a tree-cut decomposition of $G$ without increase the adhesion such that each torso of the final tree-cut decomposition of $G$ is obtained from a torso of the tree-cut decomposition of some component of $G$ by adding at most 1 isolated vertex, so the final tree-cut decomposition of $G$ satisfies the conclusion of this lemma as $H$ has no isolated vertex.

By Theorem \ref{theta-smooth}, there exists a $\theta$-smooth tree-cut decomposition $(T,\X)$ of $G$.
Note that we may assume that every $\theta$-cell in $(T,\X)$ consists of a node of $T$, since we can contract a $\theta$-cell into a node without violating the condition of being $\theta$-smooth.
So the adhesion of $(T,\X)$ is at most $\theta-1 \leq \eta$.

If there exists no $\theta$-fat $\theta$-cell in $(T,\X)$, then the torso at each node $t$ has at most $3\theta-3 \leq \xi$ edges, so we are done by choosing $Z_t=E(G_t)$ and $U_t=\emptyset$ for every $t \in V(T)$.
Hence we may assume that there exists a $\theta$-fat $\theta$-cell in $(T,\X)$. 
Let $r$ be a node of $T$ contained in a $\theta$-fat $\theta$-cell.
Recall that every $\theta$-cell in $(T,\X)$ consists of a node of $T$, so the $\theta$-fat $\theta$-cell containing $r$ consists of $r$.
Let $\E_r$ be the edge-tangle in $G$ of order $\theta$ defined by the $\theta$-cell containing $r$.
Recall that the edge-tangle in $G$ defined by a $\theta$-fat $\theta$-cell is the edge-tangle stated in Lemma \ref{fat_cell_edge-tangle}.

Since $G$ is connected, there exists a function $\iota$ satisfying the conclusions of Lemma \ref{edge-tangle_tree}.
So $V(\iota(\E_r))=\{r\}$.
Now we treat $T$ as a rooted tree rooted at $r$.

We say an edge-tangle in $G$ of order $\theta$ is {\it important} if it is defined by a $\theta$-fat $\theta$-cell.
We say an important edge-tangle $\E_1$ in $G$ of order $\theta$ is an {\it ancestor} of an important edge-tangle $\E_2$ in $G$ of order $\theta$ if some node in $V(\iota(\E_1))$ is an ancestor of some node in $V(\iota(\E_2))$; we say that $\E_2$ is a {\it descendant} of $\E_1$ if $\E_1$ is an ancestor of $\E_2$.
Note that every important edge-tangle in $G$ of order $\theta$ is an ancestor and a descendant of itself.

Since $G$ does not contain an $H$-immersion, $G$ does not contain an $H'$-immersion.
By Lemma \ref{basic_exclu_imm}, for every edge-tangle $\E$ in $G$ of order $\theta$, there exist $Z_\E \subseteq E(G)$ with $\lvert Z_\E \rvert \leq \xi_0$, $U_\E \subseteq V(G)$ with $\lvert U_\E \rvert \leq h'-1$ and a cross-free family $\D_\E \subseteq \E-Z_\E$ such that for every $v \in V(G)-U_\E$, there exists $[A,B] \in \D_\E$ of order at most $d'_{\lvert U_\E \rvert+1}-1$ with $v \in A$.
Note that we may assume that every member of $\D_\E$ has order at most $d'_{\lvert U_\E \rvert+1}-1$, for otherwise we can remove this member from $\D_\E$.
In addition, there are at most $2\lvert Z_\E \rvert \leq 2\xi_0$ members of $\D_\E$ whose order in $G-Z_\E$ and in $G$ are different, and for each such member, its order in $G$ is at most $d'_{\lvert U_\E \rvert+1}-1+\lvert Z_\E \rvert \leq d'_{\lvert U_\E \rvert+1}-1+\xi_0$.

Hence for every edge-tangle $\E$ of order $\theta$, $\D_\E \subseteq \E$ is a cross-free family in $G$ such that 
	\begin{itemize}
		\item there exists $\D'_\E \subseteq \D_\E$ with $\lvert \D'_\E \rvert \leq 2\xi_0$ such that every member of $\D'_\E$ has order at most $d'_{\lvert U_\E \rvert+1}-1+\xi_0$,
		\item every member of $\D_\E-\D'_\E$ has order at most $d'_{\lvert U_\E \rvert+1}-1$, and
		\item $\bigcup_{[A,B] \in \D_\E}A \supseteq V(G)-U_\E$. 
	\end{itemize}
For each important edge-tangle $\E$ in $G$ of order $\theta$, 
	\begin{itemize}
		\item let $\C_\E=\{\E': \E'$ is an important edge-tangle in $G$ of order $\theta$ with $\D_\E \subseteq \E'$ such that there exists no important edge-tangle $\E''$ of order $\theta$ with $\D_\E \not \subseteq \E''$ such that some path in $T$ from $V(\iota(\E))$ to $V(\iota(\E'))$ intersects $V(\iota(\E''))\}$,
		\item let $\overline{\C_\E}$ be the set of all important edge-tangles of order $\theta$ not contained in $\C_\E$, and
		\item if $\E \neq \E_r$, then let $e_\E$ be the edge of $T$ with exactly one end in $V(\iota(\E))$ such that $\iota(\E)$ and $\iota(\E_r)$ belong to different components of $T-e_\E$.
	\end{itemize}
	
\noindent{\bf Claim 1:} For every important edge-tangle $\E$ in $G$ of order $\theta$, there exists $\Se_\E$ such that $\Se_\E$ is a cross-free family and is a $(\C_\E, \overline{\C_\E})$-segregator.

\noindent{\bf Proof of Claim 1:}
For every $\E' \in \overline{\C_\E}$, by the definition of $\C_\E$, there exists an edge $e$ of $T$ with exactly one end in $\iota(\E')$ such that the component of $T-e$ containing $\iota(\E')$ is disjoint from $\iota(\E'')$ for every $\E'' \in \C_\E$, so by Lemma \ref{edge-tangle_tree}, there exists a $(\C_\E,\E')$-separator, and hence there exists a minimum $(\C_\E,\E')$-separator.
For every $\E' \in \overline{\C_\E}$, let $[A_{\E'},B_{\E'}]$ be a minimum $(\C_\E,\E')$-separator. 
Let $\F = \{[A_{\E'},B_{\E'}]: \E' \in \overline{\C_\E}\}$.
Then $\F$ is a $(\C_\E,\overline{\C_\E})$-segregator.
By Lemma \ref{uncross_segregator}, there exists $\Se_\E$ such that $\Se_\E$ is a $(\C_\E,\overline{\C_\E})$-segregator and $\Se_\E$ is a cross-free family. 
$\Box$

Let $\I$ be a collection of important edge-tangles in $G$ of order $\theta$ such that 
	\begin{itemize}
		\item $\E_r \in \I$, 
		\item for every important edge-tangle $\E$ in $G$ of order $\theta$, there exists an ancestor $\E'$ of $\E$ such that $\E' \in \I$ and $\E \in \C_{\E'}$, and
		\item for any distinct members $\E_1,\E_2$ of $\I$, if $\E_1$ is an ancestor of $\E_2$, then $\E_2 \not \in \C_{\E_1}$.
	\end{itemize}
Note that $\I$ can be easily constructed by a greedy algorithm according to a breadth-search-order of $\theta$-fat $\theta$-cells.

\noindent{\bf Claim 2:} For every $\E \in \I$ and $\E' \in \overline{\C_\E}$, there exists $[C,D] \in \Se_\E$ such that $B_{e,\iota(\E')} \subseteq C$, where $e$ is the edge of $T$ with exactly one end in $V(\iota(\E'))$ such that $\bigcup_{\E'' \in \C_\E}\iota(\E'')$ and $\iota(\E')$ belong to different components of $T-e$.

\noindent{\bf Proof of Claim 2:}
Since $\Se_\E$ is a $(\C_\E,\overline{\C_\E})$-segregator, there exists $[C,D] \in \Se_\E$ such that either $[C,D]$ is a minimum $(\C_\E,\E')$-separator, or $C' \subseteq C$ for some minimum $(\C_\E,\E')$-separator $[C',D']$.
Let $e$ be the edge of $T$ incident with $V(\iota(\E'))$ such that $\bigcup_{\E'' \in \C_\E}\iota(\E'')$ and $\iota(\E')$ belong to different components of $T-e$.
By Statement 3 of Lemma \ref{edge-tangle_tree}, if $[C,D]$ is a minimum $(\C_\E,\E')$-separator, then $B_{e,\iota(\E')} \subseteq C$; if $C' \subseteq C$ for some minimum $(\C_\E,\E')$-separator $[C',D']$, then $B_{e,\iota(\E')} \subseteq C' \subseteq C$.
Hence $B_{e,\iota(\E')} \subseteq C$.
$\Box$

\noindent{\bf Claim 3:} For every $\E \in \I$, $\Se_\E$ is a $(\D_\E,\C_\E)$-guard.

\noindent{\bf Proof of Claim 3:}
Let $[A,B] \in \Se_\E$.
Since $\Se_\E$ is a $(\C_\E,\overline{\C_\E})$-segregator by Claim 1, there exists $\E' \in \overline{\C_\E}$ such that $[A,B]$ is a minimum $(\C_\E,\E')$-separator.
Note that $\E'$ is important by the definition of $\overline{\C_\E}$.
We further assume that the distance in $T$ between $V(\iota(\E'))$ and $V(\iota(\E))$ is as small as possible.
Since $\Se_\E$ is cross-free by Claim 1, it suffices to show that $\D_\E \not \subseteq \E'$.

Suppose that $\D_\E \subseteq \E'$.
Since $\E' \not \in \C_\E$, there exists an important edge-tangle $\E''$ of order $\theta$ such that $\D_\E \not \subseteq \E''$ and $V(\iota(\E''))$ intersects the path in $T$ between $V(\iota(\E))$ and $V(\iota(\E'))$.
So $\E'' \in \overline{\C_\E}$.
Let $e''$ be the edge of $T$ incident with $V(\iota(\E''))$ such that $\iota(\E)$ and $\iota(\E'')$ belong to different components of $T-e''$.
So $\iota(\E'')$ and $\bigcup_{\E''' \in \C_\E}\iota(\E''')$ belong to different components of $T-e''$.
By Claim 2, there exists $[C,D] \in \Se_\E$ such that $B_{e'',\iota(\E'')} \subseteq C$.

Since $\E' \in \overline{\C_\E}$, by the definition of $\C_\E$, there exists an edge $e'$ of $T$ with exactly one end in $V(\iota(\E'))$ such that the component of $T-e'$ containing $\iota(\E')$ is disjoint from $V(\iota(\E'''))$ for every $\E''' \in \C_\E$.
Since $[A,B]$ is a minimum $(\C_\E,\E')$-separator, by Statement 3 of Lemma \ref{edge-tangle_tree}, $B_{e',\iota(\E')} \subseteq A$.

Note that $e''$ belongs to the path in $T$ from $e'$ to $\iota(\E)$.
So $B_{e',\iota(\E')} \subseteq B_{e'',\iota(\E'')} \subseteq C$.
Hence $B_{e',\iota(\E')} \subseteq A \cap C$.
But $\iota(\E')$ is a $\theta$-fat $\theta$-cell, $B_{e',\iota(\E')} \neq \emptyset$.
Since $\Se_\E$ is a cross-free family, $[A,B]=[C,D]$.
Since $[A,B]$ is a minimum $(\C_\E,\E')$-separator, $[A,B] \in \E'''$ for every $\E''' \in \C_\E$.
Since $[A_{e'',\iota(\E'')},B_{e'',\iota(\E'')}] \in \E''$ (by Statement 1 of Lemma \ref{edge-tangle_tree}) and $B_{e'',\iota(\E'')} \cap B = B_{e'',\iota(\E'')} \cap D \subseteq C \cap D = \emptyset$, by (E2), $[A,B] \not \in \E''$.
So $[B,A] \in \E''$ by (E1).
Hence $[A,B]$ is a $(\C_\E,\E'')$-separator.

Since $\E'$ is chosen so that the distance in $T$ between $V(\iota(\E'))$ and $V(\iota(\E))$ is minimal, $[A,B]$ is not a minimum $(\C_\E,\E'')$-separator.
Let $[A^*,B^*]$ be a minimum $(\C_\E,\E'')$-separator such that $A^*$ is maximal.
So $\lvert [A^*,B^*] \rvert < \lvert [A,B] \rvert$.
Since $[A,B]$ is a minimum $(\C_\E,\E')$-separator, $[A^*,B^*]$ is not a $(\C_\E,\E')$-separator.
Since $[A^*,B^*]$ is a $(\C_\E,\E'')$-separator, $[A^*,B^*] \in \E'$.
If $A \subseteq A^*$, then $B_{e',\iota(\E')} \subseteq A \subseteq A^*$, so $B_{e',\iota(\E')} \cap B^*=\emptyset$, but $[A_{e',\iota(\E')},B_{e',\iota(\E')}] \in \E'$ by Statement 1 of Lemma \ref{edge-tangle_tree}, contracting (E2).
Hence $A \cup A^* \supset A^*$.

If $\lvert [A^* \cup A, B^* \cap B] \rvert \leq \lvert [A^*,B^*] \rvert$, then $[A^* \cup A, B^* \cap B] \in \E'''-\E''$ and $[B^* \cap B,A^* \cup A] \in \E''-\E'''$ for every $\E''' \in \C_\E$ by Lemma \ref{easy_tangle}, so $[A^* \cup A, B^* \cap B]$ is a minimum $(\C_\E,\E'')$-separator with $A^* \cup A \supset A^*$, a contradiction.
So $\lvert [A^* \cup A, B^* \cap B] \rvert > \lvert [A^*,B^*] \rvert$.
By the submodularity, $\lvert [A^* \cap A, B^* \cup B] \rvert < \lvert [A,B] \rvert$.
So $[A^* \cap A, B^* \cup B]$ is not a $(\C_\E,\E')$-separator.
Since $[A,B]$ is a $(\C_\E,\E')$-separator, $[A^* \cap A, B^* \cup B] \in \E'''$ for every $\E''' \in \C_\E$ by Lemma \ref{easy_tangle}.
Hence $[A^* \cap A, B^* \cup B] \in \E'$.
Since $[A^*,B^*]$ is a minimum $(\C_\E,\E'')$-separator, $B_{e',\iota(\E')} \subseteq B_{e'', \iota(\E'')} \subseteq A^*$ by Statement 3 of Lemma \ref{edge-tangle_tree}.
Hence $[A_{e',\iota(\E')},B_{e',\iota(\E')}]$, $[B,A]$, $[A^* \cap A, B^* \cup B]$ are members of $\E'$ such that $B_{e', \iota(\E')} \cap A \cap (B^* \cup B) = \emptyset$, contradicting (E2).
This proves the claim.
$\Box$

For every $\E \in \I$, since $\Se_\E$ is a $(\D_\E,\C_\E)$-guard by Claim 3, by Lemma \ref{insert_guard}, there exists a cross-free family $\D_{\E,1}$ such that
	\begin{itemize}
		\item there exists $\D'_{\E,1} \subseteq \D_{\E,1}$ with $\lvert \D'_{\E,1} \rvert \leq 2\xi_0$,
		\item every edge-cut in $\D_{\E,1}-\D'_{\E,1}$ has order at most $d'_{\lvert U_\E \rvert+1}-1$, and every edge-cut in $\D'_{\E,1}$ has order at most $d'_{\lvert U_\E \rvert+1}-1+2\xi_0$,
		\item $\bigcup_{[A,B] \in \D_{\E,1}}A \supseteq V(G)-U_\E$,
		\item $\D_{\E,1} \subseteq \E'$ for every $\E' \in \C_\E$, and
		\item for every $[C,D] \in \Se_\E$, there exists $[A,B] \in \D_{\E,1}$ such that $C \subseteq A$. 
	\end{itemize}
Let $[A_{e_{\E_r},\iota(\E_r)},B_{e_{\E_r},\iota(\E_r)}] = [\emptyset, V(G)]$.
For every $\E \in \I$, 
	\begin{itemize}
		\item let $\D_{\E,2} = \{[A_{e_\E,\iota(\E)},B_{e_\E,\iota(\E)}]\} \cup \{[A \cap B_{e_\E,\iota(\E)}, B \cup A_{e_\E,\iota(\E)}]: [A,B] \in \D_{\E,1}\}$, and
		\item let $\Se_\E' = \{[A_{e_\E,\iota(\E)},B_{e_\E,\iota(\E)}]\} \cup \{[A \cap B_{e_\E,\iota(\E)}, B \cup A_{e_\E,\iota(\E)}]: [A,B] \in \Se_\E\}$. 
	\end{itemize}
Note that for every $\E \in \I$, $\Se'_\E$ is a cross-free family, and $\lvert [A_{e_\E,\iota(\E)},B_{e_\E,\iota(\E)}] \rvert < \theta$ since $e_\E$ does not belong to a $\theta$-cell.

\noindent{\bf Claim 4:} For distinct $\E,\E' \in \I$ such that $\E$ is an ancestor of $\E'$ and $\E' \in \overline{\C_\E}$, there exists $[C,D] \in \Se'_\E-\{[A_{e_\E,\iota(\E)}, B_{\E,\iota(\E)}]\}$ such that $B_{e_{\E'},\iota(\E')} \subseteq C$.

\noindent{\bf Proof of Claim 4:}
By Claim 2, there exists $[C_0,D_0] \in \Se_\E$ such that $B_{e_{\E'},\iota(\E')} \subseteq C_0$.
Since $\E'$ is a descendant of $\E$, $B_{e_{\E'},\iota(\E')} \subseteq B_{e_\E,\iota(\E)}$.
Hence $B_{e_{\E'},\iota(\E')} \subseteq C_0 \cap B_{e_\E,\iota(\E)}$.
We are done since $[C_0 \cap B_{e_\E,\iota(\E)}, D_0 \cup A_{e_\E,\iota(\E)}] \in \Se'_\E-\{[A_{e_\E,\iota(\E)}, B_{\E,\iota(\E)}]\}$.
$\Box$

By Lemma \ref{insert_one}, for every $\E \in \I$, there exists $\D'_{\E,2} \subseteq \D_{\E,2}$ with 
	\begin{itemize}
		\item $\{[A_{e_\E,\iota(\E)},B_{e_\E,\iota(\E)}]\} \cup \{[A \cap B_{e_\E,\iota(\E)}, B \cup A_{e_\E,\iota(\E)}]: [A,B] \in \D'_{\E,1}\} \subseteq \D'_{\E,2}$ and
		\item $\lvert \D'_{\E,2} \rvert \leq 2 \lvert [A_{e_\E,\iota(\E)},B_{e_\E,\iota(\E)}] \rvert+\lvert \D'_{\E,1} \rvert+1 \leq 2\theta+2\xi_0$
	\end{itemize} 
	such that
	\begin{itemize}
		\item[(D1)] every member of $\D_{\E,2}-\D'_{\E,2}$ has order at most $d'_{\lvert U_\E \rvert+1}-1$, 
		\item[(D2)] every member of $\D'_{\E,2}$ has order at most $d'_{\lvert U_\E \rvert+1}-1+2\xi_0+\theta \leq 2\theta$, 
		\item[(D3)] $\bigcup_{[A,B] \in \D_{\E,2}}A \supseteq V(G)-U_\E$, and
		\item[(D4)] for every $[C,D] \in \Se'_\E-\{A_{e_\E,\iota(\E)}, B_{\E,\iota(\E)}\}$, there exists $[A,B] \in \D_{\E,2}$ such that $C \subseteq A$. 
	\end{itemize}
For every $\E \in \I$, we define the following.
	\begin{itemize}
		\item Let $T_\E$ be the component of $T-\{e_{\E'}: \E' \in \I-\{\E_r\}\}$ containing $\iota(\E)$, and let $r_\E$ be the root of $T_\E$. 
		\item Let $G_\E = G[X_{V(T_\E)}]$. 
		\item Let $\{T_C: [C,D] \in \D_{\E,2}-\{[A_{e_\E,\iota(\E)},B_{e_\E,\iota(\E)}]\}\}$ be a set of $\lvert \D_{\E,2}-\{[A_{e_\E,\iota(\E)},B_{e_\E,\iota(\E)}]\} \rvert$ copies of $T_\E$.
		\item For each $[C,D] \in \D_{\E,2}-\{[A_{e_\E,\iota(\E)},B_{e_\E,\iota(\E)}]\}$ and $t \in V(T_\E)$, let $t_C$ be the copy of $t$ in $T_C$.
		\item Let $T^\E$ be the tree obtained from $\bigcup_{[C,D] \in \D_{\E,2}-\{[A_{e_\E,\iota(\E)},B_{e_\E,\iota(\E)}]\}}T_C$ by adding a new node $t_\E$, and for each $[C,D] \in \D_{\E,2}-\{[A_{e_\E,\iota(\E)},B_{e_\E,\iota(\E)}]\}$, adding an edge $t_\E {r_\E}_C$.
		\item Let $X^\E_{t_\E} = \bigcap_{[C,D] \in \D_{\E,2}}D$.
		\item For every $[C,D] \in \D_{\E,2}-\{[A_{e_\E,\iota(\E)},B_{e_\E,\iota(\E)}]\}$ and $t \in V(T_\E)$, let $X^\E_{t_C} = X_t \cap C$.
		\item Let $\X^\E = (X^\E_t: t \in V(T^\E))$.
	\end{itemize}
Note that $\{V(G_\E): \E \in \I\}$ is a weak partition of $V(G)$. 
(That is, members of $\{V(G_\E): \E \in \I\}$ are disjoint, and the union of all members of $\{V(G_\E): \E \in \I\}$ equals $V(G)$.)

\noindent{\bf Claim 5:} For every $\E \in \I$, $(T^\E,\X^\E)$ is a tree-cut decomposition of $G_\E$.

\noindent{\bf Proof of Claim 5:}
Clearly, $\X^\E$ consists of pairwise disjoint sets.
To show that $(T^\E,\X^\E)$ is a tree-cut decomposition of $G_\E$, it suffices to show that $\bigcup_{t \in V(T^\E)}X^\E_t = X_{V(T_\E)}$.

Note that $\bigcup_{t \in V(T^\E)}X^\E_t = X^\E_{t_\E} \cup (X_{V(T_\E)} \cap \bigcup_{[C,D] \in \D_{\E,2}-\{[A_{e_\E,\iota(\E)},B_{e_\E,\iota(\E)}]\}}C)$
$$= \allowbreak (B_{e_\E,\iota(\E)} \allowbreak \cap \allowbreak \bigcap_{[C,D] \in \D_{\E,2}-\{[A_{e_\E,\iota(\E)},B_{e_\E,\iota(\E)}]\}}D) \cup (X_{V(T_\E)} \cap \bigcup_{[C,D] \in \D_{\E,2}-\{[A_{e_\E,\iota(\E)},B_{e_\E,\iota(\E)}]\}}C).$$
Since $X_{V(T_\E)} \subseteq B_{e_\E,\iota(\E)}$, $\bigcup_{t \in V(T^\E)}X^\E_t \supseteq X_{V(T_\E)}$.

Suppose $\bigcup_{t \in V(T^\E)}X^\E_t \neq X_{V(T_\E)}$.
Then there exists $v \in \bigcup_{t \in V(T^\E)}X^\E_t - X_{V(T_\E)}$.
Note that $v \in (B_{e_\E,\iota(\E)} \cap \bigcap_{[C,D] \in \D_{\E,2}-\{[A_{e_\E,\iota(\E)},B_{e_\E,\iota(\E)}]\}}D)-X_{V(T_\E)}$.
Since $v \in B_{e_\E,\iota(\E)}-X_{V(T_\E)}$, there exists $\E' \in \I-\{\E\}$ such that $\E'$ is a descendant of $\E$ and $v \in X_{V(T_{\E'})} \subseteq B_{e_{\E'},\iota(\E')}$.
Since $\E$ and $\E'$ are distinct members of $\I$, $\E' \in \overline{\C_\E}$.
By Claim 4 and (D4), there exists $[A^*,B^*] \in \D_{\E,2}$ such that $B_{e_{\E'},\iota(\E')} \subseteq A^*$.
So $v \in A^*$.
Since $B_{e_{\E'},\iota(\E')} \subseteq A^*$ and $\E$ is an ancestor of $\E'$, $[A^*,B^*] \in \D_{\E,2}-\{[A_{e_\E,\iota(\E)},B_{e_\E,\iota(\E)}]\}$, so $v \not \in \bigcap_{[C,D] \in \D_{\E,2}-\{[A_{e_\E,\iota(\E)},B_{e_\E,\iota(\E)}]\}}D$, a contradiction.
$\Box$

For $\E_1,\E_2 \in \I$, we say that $\E_2$ is a {\it successor} of $\E_1$ if $\E_1 \neq \E_2$, $\E_2$ is a descendant of $\E_1$, and there exists no $\E' \in \I-\{\E_1,\E_2\}$ such that $\E'$ is a descendant of $\E_1$ and is an ancestor of $\E_2$.
For each $\E \in \I$ and successor $\E'$ of $\E$, by Claim 4 and (D4), there exists uniquely $[C_{\E,\E'},D_{\E,\E'}] \in \D_{\E,2}-\{[A_{e_\E,\iota(\E)},B_{e_\E,\iota(\E)}]\}$ such that $B_{e_{\E'},\iota(\E')} \subseteq C_{\E,\E'}$.
(Note that the uniqueness follows from the fact that $\D_{\E,2}$ is cross-free.)

For each $\E \in \I-\{\E_r\}$, let $r'_\E$ be the end of $e_\E$ other than $r_\E$.

Finally, we define the following.
	\begin{itemize}
		\item Define $T^*$ to be the tree obtained from $\bigcup_{\E \in \I}T^\E$ by for each $\E \in \I$ and each successor $\E' \in \I$ of $\E$, adding an edge ${r'_{\E'}}_{C_{\E,\E'}}t_{\E'}$. 
		\item For each $\E \in \I$ and $t \in V(T^\E)$, let $X^*_t=X^\E_t$.
		\item Define $\X^*=(X^*_t: t \in V(T^*))$.
	\end{itemize}
Note that by Claim 5, for each $\E \in \I$, $X^*_{V(T^\E)}=X^\E_{V(T^\E)}=X_{V(T_\E)}$.
Since $\{V(G_\E): \E \in \I\}$ is a weak partition of $V(G)$, $(T^*,\X^*)$ is a tree-cut decomposition of $G$.
We shall prove that $(T^*,\X^*)$ satisfies the conclusion of this lemma.

\noindent{\bf Claim 6:} The adhesion of $(T^*,\X^*)$ is at most $\eta$.

\noindent{\bf Proof of Claim 6:}
Let $e \in E(T^*)$.
If $e$ is not an edge of $\bigcup_{\E \in \I}T^\E$, then $e={r'_{\E'}}_{C_{\E,\E'}}t_{\E'}$ for some $\E,\E' \in \I$, where $\E'$ is a successor of $\E$, so $\adh_{(T^*,\X^*)}(e) = \adh_{(T,\X)}(r_{\E'}r'_{\E'})$ by Claim 5, and hence $\lvert \adh_{(T^*,\X^*)}(e) \rvert = \lvert \adh_{(T,\X)}(r_{\E'}r'_{\E'}) \rvert < \theta \leq \eta$.
If $e$ is an edge of $\bigcup_{\E' \in \I}T^{\E'}$ incident with $t_\E$ for some $\E \in \I$, then by the definition of $C_{\E,\E'}$ and $\D_{\E,2}-\{[A_{e_\E,\iota(\E)},B_{e_\E,\iota(\E)}]\}$, we know that $\adh_{(T^*,\E^*)}(e)$ is the set of edges between $C$ and $D$ for some $[C,D] \in \D_{\E,2}-\{[A_{e_\E,\iota(\E)},B_{e_\E,\iota(\E)}]\}$, so it has size at most $\eta$ by (D1) and (D2).

So we may assume that $e$ is an edge of $T^\E$ for some $\E \in \I$ not incident with $t_\E$.
So there exists an edge $e'$ of $T$ such that $e$ is a copy of $e'$.
By Claim 5 and the definition of $\C_{\E,\E'}$ for successors $\E'$ of $\E$, $\adh_{(T^*,\X^*)}(e)$ is the set of edges between $A_{e',r} \cap C$ and $B_{e',r} \cup D$ for some $[C,D] \in \D_{\E,2}-\{[A_{e_\E,\iota(\E)},B_{e_\E,\iota(\E)}]\}$, so it has size at most $\lvert [A_{e',r},B_{e',r}] \rvert + \lvert [C,D] \rvert \leq \lvert [A_{e',r},B_{e',r}] \rvert+2\theta$ by (D1) and (D2).
Recall that we assume that every $\theta$-cell in $(T,\X)$ is a node, so $\lvert [A_{e',r},B_{e',r}] \rvert \leq \theta$.
Hence $\lvert \adh_{(T^*,\X^*)}(e) \rvert \leq \lvert [A_{e',r},B_{e',r}] \rvert+2\theta \leq 3\theta \leq \eta$.
$\Box$

\noindent{\bf Claim 7:} For every $t \in V(T^*)-\{t_\E: \E \in \I\}$, the torso at $t$ has at most $\xi$ edges.

\noindent{\bf Proof of Claim 7:}
Let $t \in V(T^*)-\{t_\E: \E \in \I\}$.
So $t \in V(T^\E)-\{t_\E\}$ for some $\E \in \I$.
Hence there exist $t' \in V(T)$ and $[C,D] \in \D_{\E,2}-\{[A_{e_\E,\iota(\E)},B_{e_\E,\iota(\E)}]\}$ such that $t=t'_C$.
Note that every edge of the torso at $t$ in $(T^*,\X^*)$ is either an edge of the torso at $t'$ in $(T,X)$ incident with $C$ or an edge of $G$ between $C$ and $D$.
If $t'$ is not contained in any $\theta$-fat $\theta$-cell, then the torso at $t$ in $(T^*,\X^*)$ has at most $3\theta-3+2\theta \leq \xi$ edges by (D1) and (D2).

So we may assume that $t'$ is contained in some $\theta$-fat $\theta$-cell $L$.
Recall that we assume that every $\theta$-cell consists of a node, so $V(L)=\{t'\}$.
Hence there exists an edge-tangle $\E_L$ defined by $L$ as mentioned in Lemma \ref{fat_cell_edge-tangle}.
Note that $\iota(\E_L) \subseteq T_\E$.
So $\E_L \in \C_\E$ and $\D_{\E,1} \subseteq \E_L$.
Since $[C,D] \in \D_{\E,2}-\{[A_{e_\E,\iota(\E)},B_{e_\E,\iota(\E)}]\}$, $[C,D]=[C' \cap B_{e_\E,\iota(\E)}, D' \cup A_{e_\E,\iota(\E)}]$ for some $[C',D'] \in \D_{\E,1}$.
Since $[C',D'] \in \D_{\E,1} \subseteq \E_L$, by Lemma \ref{fat_cell_edge-tangle}, $C'$ is incident with at most $\lvert [C',D'] \rvert$ edges of the torso at $\iota(\E_L)$ in $(T,\X)$.
Hence $C=C' \cap B_{e_\E,\iota(\E)}$ is incident with at most $\lvert [C',D'] \rvert \leq \theta$ edges of the torso at $\iota(\E_L)=L$ in $(T,\X)$.
Since every edge of the torso at $t$ in $(T^*,\X^*)$ is either an edge of the torso at $t'$ in $(T,\X)$ incident with $C$ or an edge of $G$ between $C$ and $D$, the torso at $t$ has at most $\theta + \lvert [C,D] \rvert \leq 3\theta \leq \xi$ edges by (D1) and (D2).
$\Box$

Since $H$ has no isolated vertex, Claim 7 implies that this lemma holds for every $t \in V(T^*)-\{t_\E: \E \in \I\}$, as one can define $Z_t$ to be the set of all edges in the torso at $t$ and define $U_t=\emptyset$.

\noindent{\bf Claim 8:} For every $\E \in \I$, if $G_{t_\E}$ is the torso at $t_\E$, then there exists $Z \subseteq E(G)$ with $\lvert Z \rvert \leq \xi$ and $U \subseteq X_{t_\E}$ with $\lvert U \rvert \leq h'-1$ such that for every $v \in V(G_{t_\E})-U$, the degree of $v$ in $G_{t_\E}-Z$ is at most $d'_{\lvert U \rvert+1}-1$.

\noindent{\bf Proof of Claim 8:}
Let $\E \in \I$.
Let $G_{t_\E}$ be the torso at $t_\E$.
Define $U=U_\E \cap X_{t_\E}$.
Define $Z=\{e \in E(G): e$ is between $C$ and $D$ for some $[C,D] \in \D'_{\E,2}\}$.
Note that $\lvert U \rvert \leq \lvert U_\E \rvert \leq h'-1$.
By (D2), $\lvert Z \rvert \leq \lvert \D'_{\E,2} \rvert \cdot 2\theta \leq (2\theta+2\xi_0) \cdot 2\theta \leq 8\theta^2=\xi$.

Note that $G_{t_\E}$ is obtained from $G$ by for each $[C,D] \in \D_{\E,2}$, identifying $C$ into a vertex and deleting resulting loops.
By the definition of $Z$, every peripheral vertex corresponding to a member of $\D'_{\E,2}$ has degree 0 in $G_{t_\E}-Z$.
Note that $0 \leq d'_{\lvert U \rvert+1}-1$ since $H$ has no isolated vertex.
By (D1), every peripheral vertex corresponding to a member of $\D_{\E,2}-\D_{\E,2}'$ has degree at most $d'_{\lvert U_\E \rvert+1}-1 \leq d'_{\lvert U \rvert+1}-1$.
By (D3), every non-peripheral vertex of $G_{t_\E}$ is in $X_{t_\E} \cap \bigcap_{[C,D] \in \D_{\E,2}}D \subseteq U$.
This proves the claim.
$\Box$

Claim 8 completes the proof of this lemma.
\end{pf}

\bigskip

Now we are ready to prove a global decomposition theorem for excluding $H$-immersions in graphs with no edge-cut of order 3, where $H$ is allowed to have isolated vertices.

\begin{theorem} \label{global}
For any positive integers $d,h$, there exist integers $\eta=\eta(d,h)$ and $\xi=\xi(d,h)$ such that the following holds.
Let $H$ be a graph on $h$ vertices with maximum degree $d$. 
Let $G$ be a graph with no edge-cut of order exactly 3 such that $G$ does not contain an $H$-immersion.
Define $H'=H$ if $H$ is non-exceptional; otherwise, define $H'$ to be a graph obtained from $H$ by subdividing one edge.
Then there exists a tree-cut decomposition $(T,\X)$ of $G$ of adhesion at most $\eta$ such that for every $t \in V(T)$, there exists $Z_t \subseteq E(G)$ with $\lvert Z_t \rvert \leq \xi$ such that if $G_t$ is the torso at $t$, then there exists a nonnegative integer $k_t \leq d$ such that 
	\begin{enumerate}
		\item the number of vertices of degree at least $k_t$ in $G_t-Z_t$ is less than the number of vertices of degree at least $k_t$ in $H'$, 
		\item every vertex of $G_t$ of degree at least $k_t$ in $G_t-Z_t$ is a non-peripheral vertex of $G_t$,
		\item if $\lvert V(T) \rvert=1$ or $t$ is not a leaf, then every vertex in $X_t$ has degree at least $k_t$ in $G_t-Z_t$, and
		\item if $t$ is a leaf and $\lvert V(T) \rvert \geq 2$, then $\lvert X_t \rvert \leq 1$.
	\end{enumerate}
\end{theorem}

\begin{pf}
Let $\eta_1=\max_{1 \leq i \leq h+1}\eta_{\ref{global_no_iso_5}}(d,i)$ and $\xi_1=\max_{1 \leq i \leq h+1}\xi_{\ref{global_no_iso_5}}(d,i)$, where $\eta_{\ref{global_no_iso_5}}$ and $\xi_{\ref{global_no_iso_5}}$ are the integers $\eta$ and $\xi$ mentioned in Lemma \ref{global_no_iso_5}, respectively.
Define $\eta = \eta_1 + \xi_1$ and $\xi = \xi_1$.
Note that $\eta \geq d+\xi_1$ by Lemma \ref{global_no_iso_5}.

\noindent{\bf Claim 1:} There exists a tree-cut decomposition $(T,\X)$ of $G$ of adhesion at most $\eta_1$ such that for every $t \in V(T)$, there exists $Z_t \subseteq E(G)$ with $\lvert Z_t \rvert \leq \xi_1$ and such that if $G_t$ is the torso at $t$, then there exists a nonnegative integer $k_t \leq d$ such that 
	\begin{itemize}
		\item the number of vertices of degree at least $k_t$ in $G_t-Z_t$ is less than the number of vertices of degree at least $k_t$ in $H'$, and
		\item every vertex of $G_t$ of degree at least $k_t$ in $G_t-Z_t$ is a non-peripheral vertex of $G_t$.
	\end{itemize}

\noindent{\bf Proof of Claim 1:}
We first assume that $\lvert V(G) \rvert < \lvert V(H') \rvert$.
Define $(T,\X)$ to be the tree-cut decomposition such that $T$ consists of a single node $t$ and define $X_t=V(G)$ and $\X=(X_t)$.
Define $Z=\emptyset$.
So the torso at the unique node is $G$, and there exists no peripheral vertex of the torso.
Note that the number of vertices of $G-Z$ with degree at least zero is $\lvert V(G) \rvert$, and the number of vertices of $H'$ with degree at least zero is $\lvert V(H') \rvert$.
We are done by choosing $k_t=0$ since $\lvert V(G) \rvert <\lvert V(H') \rvert$.

Hence we may assume that $\lvert V(G) \rvert \geq \lvert V(H') \rvert$.
Since $G$ does not contain an $H'$-immersion, $H'$ contains at least one edge.
Let $H^*$ be the subgraph of $H'$ induced by all non-isolated vertices of $H'$.

Since $H^*$ contains all non-isolated vertices of $H'$ and $\lvert V(G) \rvert \geq \lvert V(H') \rvert$, $G$ does not contain an $H^*$-immersion, for otherwise $G$ contains an $H'$-immersion.
Let $(d_1^*,d_2^*,...,d^*_{h^*})$ be the degree sequence of $H^*$, where $h^*=\lvert V(H^*) \rvert$.
Note that $d_1^*=d$ and $1 \leq h^* \leq \lvert V(H') \rvert \leq h+1$.
By Lemma \ref{global_no_iso_5}, there exists a tree-cut decomposition $(T,\X)$ of $G$ of adhesion at most $\eta_1$ such that for every $t \in V(T)$, there exist $Z_t \subseteq E(G)$ and $U_t \subseteq X_t$ with $\lvert Z_t \rvert \leq \xi_1$ and $\lvert U_t \rvert \leq h^*-1$ such that if $G_t$ is the torso at $t$, then for every $v \in V(G_t-Z_t)-U_t$, the degree of $v$ in $G_t-Z_t$ is at most $d^*_{\lvert U_t \rvert+1}-1$. 
Since for every $t \in V(T)$, $U_t \subseteq X_t$, every peripheral vertex of $G_t$ has degree less than $d^*_{\lvert U_t \rvert+1}$ in $G_t-Z_t$, and the number of vertices of $G_t-Z_t$ with degree at least $d^*_{\lvert U_t \rvert+1}$ is at most $\lvert U_t \rvert$, but the number of vertices of $H^*$ with degree at least $d^*_{\lvert U_t \rvert+1}$ is at least $\lvert U_t \rvert+1$.
Since $H^*$ is a subgraph of $H'$, the claim follows by defining $k_t=d^*_{\lvert U_t \rvert+1}$.
$\Box$

For every node $t \in V(T)$, let $S_t$ be the set of vertices in $X_t$ of degree in $G_t-Z_t$ less than $k_t$.
Define $T'$ to be the tree obtained from $T$ by for every node $t \in V(T)$, attaching $\lvert S_t \rvert$ leaves to $t$.
Note that there exists an obvious bijection between $S_t$ and the new leaves attached to $t$.
For every $t \in V(T)$, define $X'_t = X_t-S_t$; for every $t \in V(T)$ and a leaf $t'$ in $V(T')-V(T)$ adjacent to $t$, define $X'_{t'}$ to be the set consisting of the vertex in $S_t$ corresponding to $t'$.
Define $\X' = (X'_t: t \in V(T'))$.
For each $t' \in V(T')-V(T)$, define $Z_{t'}=Z_t$ and define $k_{t'}=k_t$, where $t$ is the node of $T$ adjacent in $T'$ to $t'$.

Since $k_t \leq d \leq \eta-\xi_1$ for each $t \in V(T)$, the adhesion of $(T',\X')$ is at most $\eta$.
For every $t \in V(T')$, let $G_t'$ be the torso at $t$ in $(T',\X')$.
For every $t \in V(T)$, $G'_t$ is obtained from $G_t$ by deleting loops incident with $S_t$.
For every $t' \in V(T')-V(T)$, $G'_t$ consists of at most two vertices, where each vertex has degree at most $k_{t'}-1 = k_t-1$ after deleting $Z_{t'}$, where $t$ is the neighbor of $t'$ in $T'$.
So $(T',\X')$ satisfies Statements 1-3 of this theorem.

Hence we may assume that $\lvert V(T') \rvert \geq 2$, for otherwise we are done.
Since $(T',\X')$ satisfies Statement 2 and $\lvert V(T') \rvert \geq 2$, we know for every $t \in V(T')$, $G'_t$ contains a peripheral vertex of degree at least 0, so $k_t \geq 1$.

For every leaf $t$ of $T'$ with $\lvert X'_t \rvert \geq 2$, $t$ is a leaf of $T$ with $S_t = \emptyset$, so $X_t=X_t'$ and every vertex in $X_t'$ has degree at least $k_t$ in $G_t-Z_t = G'_t-Z_t$.
Let $T^*$ be the tree obtained from $T'$ by for each leaf $t$ of $T'$ with $\lvert X'_t \rvert \geq 2$, attaching a leaf adjacent to $t$.
For each $t \in V(T^*)-V(T')$, define $X^*_t = \emptyset$; for each $t \in V(T')$, define $X^*_t=X'_t$.
Let $\X^* = (X^*_t: t \in V(T^*))$.
Then $(T^*,\X^*)$ is a tree-cut decomposition of $G$ with the same adhesion as $(T',\X')$.

For each $t \in V(T^*)-V(T')$, let $Z_t = \emptyset$ and $k_t=d$.
For each $t \in V(T^*)$, let $G^*_t$ be the torso at $t$.
For each $t \in V(T^*)-V(T')$, $G_t^*$ consists of an isolated vertex, and this vertex is a peripheral vertex.
Note that if $t \in V(T')$ with $G^*_t \neq G'_t$, then $t$ is a leaf of $T'$ with $\lvert X'_t \rvert \geq 2$ and $S_t = \emptyset$, and $G^*_t$ is obtained from $G_t'$ by adding an isolated vertex.
Recall that $k_t \geq 1$ for every $t \in V(T')$, so $(T^*,\X^*)$ satisfies Statements 1 and 2.
If $t \in V(T^*)$ is not a leaf in $T^*$ but is a leaf in $T'$, then $t \in V(T')$, $\lvert X'_t \rvert \geq 2$, and $S_t = \emptyset$.
Since $(T',\X')$ satisfies Statement 3, so does $(T^*,\X^*)$.
Since every leaf $t$ of $T^*$ either is a leaf of $T'$ with $\lvert X'_{t'} \rvert=1$ or satisfies $X^*_t = \emptyset$, $(T^*,\X^*)$ satisfies Statement 4.
This proves the theorem.
\end{pf}

\section{Simple byproducts} \label{sec:app}

Let $G$ be a graph.
Let $(T,\X=(X_t: t\in V(T)))$ be a tree-cut decomposition of $G$.
We say that $(T,\X)$ is a {\it weak carving} if $\lvert X_t \rvert=1$ for every leaf $t$ of $T$, and $X_t = \emptyset$ for every non-leaf node $t$ of $T$.
We say that $(T,\X)$ is a {\it carving} if it is a weak carving such that every node of $T$ is of degree 1 or 3.
The {\it carving-width} of $G$ is the minimum $w$ such that $G$ admits a tree-cut decomposition that is a carving of adhesion at most $w$.

The first objective of this section is to prove Proposition \ref{width_relation_intro}.

\begin{lemma} \label{torso_width_carving}
Let $G$ be a graph.
Let $w$ be a nonnegative integer.
	\begin{enumerate}
		\item If the tree-cut torso-width of $G$ is at most $w$, then the carving-width of $G$ is at most $w$.
		\item If the carving-width of $G$ is at most $w$ and $G$ is loopless, then the tree-cut torso-width of $G$ is at most $3w/2$.
	\end{enumerate}
\end{lemma}

\begin{pf}
We first assume that the tree-cut torso-width of $G$ is at most $w$.
So there exists a tree-cut decomposition $(T,\X=(X_t: t \in V(T)))$ of $G$ of torso-width at most $w$.
Let $T'$ be the tree obtained from $T$ by for each $t \in V(T)$, attaching $\lvert X_t \rvert$ leaves incident with $t$.
Note that there exists an bijection between $X_t$ and the set of new leaves attached to $t$.
For every $t' \in V(T')-V(T)$, $t'$ is a leaf attached to a node $t \in V(T)$, and we let $X'_{t'}$ be the set consisting of the vertex in $X_t$ corresponding to $t'$; for every $t' \in V(T)$, let $X'_{t'} = \emptyset$.
Let $\X' = (X'_t: t\in V(T))$.
Then $(T',\X')$ is a tree-cut decomposition of $G$ such that for every $t' \in V(T')$, the torso at $t'$ in $(T',\X')$ is obtained from the torso at $t$ in $(T,\X)$ for some $t \in V(T)$ by identifying vertices and deleting some loops. 
So the torso-width of $(T',\X')$ is at most $w$ and $(T',\X')$ is a weak carving.
We shall modify $T'$ to be a cubic tree without increasing the torso-width.
Let $t$ be a node of $T'$ such that the degree of $t$ in $T'$ is at least four.
Let $x,y$ be neighbors of $t$ and let $T_x,T_y$ be the components of $T-t$ containing $x,y$ respectively.
Let $T''$ be the tree obtained from $T$ by deleting the edges $tx,ty$ and adding a new vertex $t^*$ and new edges $t^*x,t^*y,t^*t$.
For every $t \in V(T')$, let $X''_t=X_t$; let $X''_{t^*} = \emptyset$.
Let $\X'' = (X''_t: t \in V(T''))$.
Then $(T'',\X'')$ is a weak carving with torso-width at most the torso-width of $(T',\X')$ such that the number of edges of $T''$ incident with nodes of degree greater than 3 is smaller than the number of edges of $T'$ incident with nodes of degree greater than 3.
So by repeatedly applying this process, we may assume that every node of $T'$ has degree at most three.
Since $(T',\X')$ is a weak carving, if there exists $t \in V(T')$ of degree 2 in $T'$, then we can contract an edge of $T'$ incident with $t$ to obtain a weak carving with fewer nodes of degree 2 without increasing the torso-width.
So by repeatedly applying this process, we obtain a carving of torso-width at most $w$ and hence of adhesion at most $w$.

Now we assume that $G$ is loopless and the carving-width of $G$ is at most $w$.
Let $(T^*,\X^*)$ be a carving with adhesion at most $w$.
For every $t \in V(T^*)$, let $G_t$ be the torso at $t$ in $(T^*,\X^*)$.
If $t$ is a leaf, then $G_t$ has at most $w$ edges since $G$ is loopless.
If $t$ is not a leaf, then for each edge of $G_t$, there exist two components $T_1,T_2$ of $T-t$ such that one end is in $X_{V(T_1)}$ and one end is in $X_{V(T_2)}$; since the adhesion of $(T^*,\X^*)$ is at most $w$ and $t$ has degree three in $T$, there are at most $3w/2$ edges in $G_t$.
So the torso-width of $(T^*,\X^*)$ is at most $3w/2$.
\end{pf}

\begin{lemma} \label{min_torso_width}
Let $G$ be a graph of tree-cut torso-width $w$.
Then for every integer $\theta$ with $w+1 \leq \theta \leq \lvert E(G) \rvert$, $G$ admits a $\theta$-smooth tree-cut decomposition $(T,\X)$ of torso-width $w$ with no $k$-cell for every $w+1 \leq k \leq \theta$.
\end{lemma}

\begin{pf}
Let $(T,\X)$ be a tree-cut decomposition of $G$ of lexicographically minimum $\theta$-signature.
By Theorem \ref{theta-smooth}, $(T,\X)$ is $\theta$-smooth.

Suppose there exists a $k$-cell in $(T,\X)$ for some $k$ with $w+1 \leq k \leq \theta$.
Since the tree-cut torso-width of $G$ is $w$, there exists a tree-cut decomposition $(T',\X')$ of $G$ such that no torso at a node of $T'$ in $(T',\X')$ has at least $w+1$ edges.
So $(T',\X')$ has adhesion at most $w$ and has no $r$-cell for every $w+1 \leq r \leq \lvert E(G) \rvert$.
Hence the $\theta$-signature of $(T',\X')$ is smaller than the $\theta$-signature of $(T,\X)$, a contradiction.

So there exists no $k$-cell in $(T,\X)$ for every $w+1 \leq k \leq \theta$.
If there exists $t \in V(T)$ such that the torso at $t$ has more than $w$ edges, then $t$ is contained in a $(w+1)$-cell, a contradiction.
So the torso-width of $(T,\X)$ is at most $w$.
\end{pf}

\begin{lemma} \label{edge-tangle_duality_1}
Let $G$ be a graph.
Let $w$ be a positive integer.
If the tree-cut torso-width of $G$ is at most $w$, then there exists no edge-tangle of order at least $w+1$ in $G$. 
\end{lemma}

\begin{pf}
If there exists an edge-tangle $\E$ of order $\theta$ in $G$ for some $\theta \geq w+1$, then for every tree-cut decomposition $(T,\X)$ of $G$, by Lemma \ref{weak_tangle_cell}, there exists a $\theta$-cell in $(T,\X)$, so the tree-cut torso-width of $G$ is at least $\theta \geq w+1$ by Lemma \ref{min_torso_width}, a contradiction.
\end{pf}

\begin{lemma} \label{edge-tangle_duality_2}
Let $G$ be a graph.
Let $\theta$ be a positive integer.
If there exists no edge-tangle of order $\theta$ in $G$, then the tree-cut torso-width is at most $3\theta-3$. 
\end{lemma}

\begin{pf}
Let $w$ the tree-cut torso with of $G$.
Let $(T,\X)$ be a tree-cut decomposition of $G$ with lexicographically minimum $\theta$-signature.
By Theorem \ref{theta-smooth}, $(T,\X)$ is $\theta$-smooth.
If $w \geq 3\theta-2$, then there exists a node $t$ of $T$ such that the torso at $t$ in $(T,\X)$ has at least $3\theta-2$ edges, so $t$ is contained in a $\theta$-cell whose torso has at least $3\theta-2$ edges, and hence there exists an edge-tangle of order $\theta$ by Lemma \ref{fat_cell_edge-tangle}, a contradiction.
Hence $w \leq 3\theta-3$.
\end{pf}

\bigskip

\noindent {\bf Proof of Proposition \ref{width_relation_intro}:}
It immediately follows from Lemmas \ref{torso_width_carving}, \ref{edge-tangle_duality_1} and \ref{edge-tangle_duality_2}.
$\Box$

\bigskip

Now we prove Corollary \ref{mw_thm_intro}.
For an integer $k$, we say that a graph $G$ is {\it $k$-simple} if every vertex of $G$ is incident with at most $k$ loops and for every pair of distinct vertices of $G$, there exist at most $k$ parallel edges between them.
The following corollary immediate implies Corollary \ref{mw_thm_intro}.

\begin{corollary} \label{mw_equiv}
For any positive integers $d,h$ with $d \geq 4$, there exists a positive integer $\xi=\xi(d,h)$ such that the following hold.
	If $G$ is an $H$-immersion free graph for some graph $H$ on at most $h$ vertices with maximum degree at most $d$, and every component of $G$ does not have an edge-cut of order $k$ for every $3 \leq k \leq d-1$, then there exists a tree-cut decomposition $(T,\X)$ of $G$ of adhesion at most $\xi$ such that 
	\begin{enumerate}
		\item the maximum degree of $T-S$ is at most $\xi$, where $S=\{e \in E(T): \lvert \adh_{(T,\X)}(e) \rvert \leq 2\}$, and 
		\item $\lvert X_t \rvert \leq \xi$ for every $t \in V(T)$.
	\end{enumerate}
In particular, 
	\begin{enumerate}
		\item $G$ has tree-cut width at most $2\xi$, and
		\item if every component of $G$ is $d$-edge-connected, and there exists an integer $k$ such that $G$ is $k$-simple, then $G$ has tree-cut torso-width at most $(k+1)\xi^2$, carving-width at most $(k+1)\xi^2$ and maximum degree at most $(k+1)\xi^2$. 
	\end{enumerate}
\end{corollary}

\begin{pf}
Define $\xi = \max\{h,\eta_{\ref{global}}(i,j),2\xi_{\ref{global}}(i,j): i \in [d], j \in [h+1]\}$, where $\eta_{\ref{global}}$ and $\xi_{\ref{global}}$ are the integers $\eta$ and $\xi$ mentioned in Theorem \ref{global}, respectively.

Let $H$ be a graph on at most $h$ vertices with maximum degree at most $d$.
Let $G$ be an $H$-immersion free graph whose every component does not contain any edge-cut of order $k$ for every $3 \leq k \leq d-1$. 
Let $H'$ be the graph obtained from $H$ as defined in Theorem \ref{global}.
By Theorem \ref{global}, there exists a tree-cut decomposition $(T,\X)$ of $G$ of adhesion at most $\xi$ such that for every $t \in V(T)$, there exists $Z_t \subseteq E(G)$ with $\lvert Z_t \rvert \leq \xi/2$ satisfying the conclusion of Theorem \ref{global}.
In particular, for every $t \in V(T)$, every peripheral vertex of the torso at $t$ has degree less than $d$ after deleting $Z_t$.
Since every component of $G$ does not contain an edge-cut of order $k$ for every $3 \leq k \leq d-1$, we know for every $t \in V(T)$, every peripheral vertex of the torso at $t$ with degree at least 3 is incident with an edge in $Z_t$, so there are at most $2\lvert Z_t \rvert \leq \xi$ peripheral vertices with degree at least 3.
Let $S = \{e \in E(T): \lvert \adh_{(T,\X)}(e) \rvert \leq 2\}$.
Hence the maximum degree of $T-S$ is at most $\xi$.
In addition, for every $t \in V(T)$, $\lvert X_t \rvert \leq \lvert V(H') \rvert-1 \leq h \leq \xi$.
	
Hence, for every $t \in V(T)$, the torso at $t$ has at most $\xi$ peripheral vertices with degree at least 3, so the 3-center of the torso at $t$ has at most $\lvert X_t \rvert + \xi \leq 2\xi$ vertices.
Since the adhesion of $(T,\X)$ is at most $\xi$, the tree-cut width of $G$ is at most $2\xi$.

When $G$ is $d$-edge-connected, $S= \emptyset$.
Hence if $G$ is $d$-edge-connected and $k$-simple for some positive integer $k$, then for every $t \in V(T)$, the torso at $t$ has at most $\xi$ peripheral vertices and at most $k \lvert X_t \rvert + k{\lvert X_t \rvert \choose 2} + \xi^2 \leq (k+1)\xi^2$ edges, so the tree-cut torso-width of $G$ is at most $(k+1)\xi^2$, and hence the carving-width of $G$ is at most $(k+1)\xi^2$ by Lemma \ref{torso_width_carving}, and hence the maximum degree of $G$ is at most $(k+1)\xi^2$.
Therefore, if every component of $G$ is $d$-edge-connected and $k$-simple, then the tree-cut torso-width, carving-width and maximum degree are at most $(k+1)\xi^2$.
\end{pf}

\bigskip

The rest of the section is dedicated to a proof of Theorem \ref{counting_basic_intro}.

Let $\xi$ be a positive integer.
A tree-cut decomposition $(T,\X)$ of a graph is {\it $\xi$-nice} if it has adhesion at most $\xi$, and for every $t \in V(T)$, the torso at $t$ is obtained from a graph on at most $\xi$ vertices by attaching leaves adjacent to $X_t$.
A graph is {\it $\xi$-nice} if it admits a $\xi$-nice tree-cut decomposition.

\begin{lemma} \label{nice_tree_cut_sep}
Let $\xi$ be a positive integer.
Let $(T,\X=(X_t:t \in V(T)))$ be a $\xi$-nice tree-cut decomposition of a graph $G$.
Then either
	\begin{enumerate}
		\item there exists an edge $e=t_1t_2$ of $T$ such that for each $i \in [2]$, there exist at least $\frac{1}{3}\lvert E(G) \rvert$ edges of $G$ incident with $A_{e,t_i}$ and $G[A_{e,t_i}]$ is a $\xi$-nice graph with at least $\frac{1}{3}\lvert E(G) \rvert-\xi$ edges, or
		\item there exist $t^* \in V(T)$ such that either $\lvert E(G[X_{t^*}]) \rvert \geq \frac{1}{9}\lvert E(G) \rvert$, or there exists a partition $\{U_1,U_2\}$ of the set of the components of $T-t^*$ such that for each $i \in [2]$, the number of edges of $G$ incident with $\bigcup_{C \in U_i}X_{V(C)}$ is at least $\frac{2}{9}\lvert E(G) \rvert$.
	\end{enumerate}
\end{lemma}

\begin{pf}
Since $(T,\X)$ is $\xi$-nice, the adhesion of $(T,\X)$ is at most $\xi$.
So we may assume that for every edge $e$ of $T$, there exists an end $t$ of $e$ such that there are less than $\frac{1}{3}\lvert E(G) \rvert$ edges of $G$ incident with $A_{e,t}$, for otherwise Statement 1 holds.
So for every edge $e$ of $T$, there uniquely exists an end $t$ of $e$ such that there are less than $\frac{1}{3}\lvert E(G) \rvert$ edges of $G$ incident with $A_{e,t}$, and we assign a direction of $e$ such that $t$ is the head of $e$.
Hence we obtain an orientation of $E(T)$.
So there exists a node $t^*$ of $T$ with out-degree 0.

We may assume that $\lvert E(G[X_{t^*}]) \rvert < \frac{1}{9}\lvert E(G) \rvert$, for otherwise we are done since Statement 2 holds.
So $T$ contains at least two nodes, and there exist at least $\frac{8}{9}\lvert E(G) \rvert$ edges incident with $V(G)-X_{t^*}$.
Let $T_1,T_2,...,T_k$ (for some positive integer $k$) be the components of $T-t^*$. 
For each $i \in [k]$, let $t_i$ be the neighbor of $t^*$ contained in $T_i$. 
Since there exist at least $\frac{8}{9}\lvert E(G) \rvert$ edges incident with $V(G)-X_{t^*}$, there exists the minimum $m$ such that there are at least $\frac{1}{3}\lvert E(G) \rvert$ edges of $G$ incident with $\bigcup_{j=1}^{m}B_{t_it^*,t_i}$.
Let $U_1 = \{T_i: i \in [m]\}$ and let $U_2 = \{T_i: m+1 \leq i \leq k\}$.
Since for each $i \in [k]$, there are less than $\frac{1}{3}\lvert E(G) \rvert$ edges incident with $B_{t_it^*,t_i}$, by the minimality of $m$, there are at least $\frac{1}{3}\lvert E(G) \rvert$ edges of $G$ incident with $\bigcup_{C \in U_1}X_C$ and there are at most $\frac{2}{3}\lvert E(G) \rvert$ edges of $G$ incident with $\bigcup_{C \in U_1}X_{V(C)}$.
Hence there are at least $(\frac{8}{9}-\frac{2}{3})\lvert E(G) \rvert = \frac{2}{9}\lvert E(G) \rvert$ edges of $G$ incident with $\bigcup_{C \in U_2}X_{V(C)}$.
\end{pf}

\begin{lemma} \label{counting_nice}
For every positive integer $\xi$, there exists an integer $c=c(\xi)$ such that for every positive integer $m$, the number of unlabelled $m$-edge $\xi$-nice graphs with no isolated vertex is at most $c^m/m^{3\xi^2+4}$.
\end{lemma}

\begin{pf}
Let $a=10^{16\xi^2}$. 
Define $c=a^{4a}$. 

We shall prove this lemma by induction on $m$.
Note that every $m$-edge graph with no isolated vertex has at most $2m$ vertices.
So there are at most ${{2m \choose 2}+2m-1 \choose m} \leq (3m^2)^m$ vertex-labelled (and hence unlabelled) $m$-edge graphs with no isolated vertex.

When $m \leq a$, there exist at most $(3m^2)^m \leq (3a^2)^a \leq \frac{c}{a^{3\xi^2+4}} \leq \frac{c^m}{m^{3\xi^2+4}}$ such graphs, so we are done. 
Hence we may assume that $m > a$ and for every positive integer $x$ with $x<m$, the number of unlabelled $x$-edge $\xi$-nice graphs with no isolated vertex is at most $c^x/x^{3\xi^2+4}$.

By Lemma \ref{nice_tree_cut_sep}, every $m$-edge $\xi$-nice graph $G$ with no isolated vertex admits a $\xi$-nice tree-cut decomposition $(T,\X)$ satisfying one of the following.
	\begin{itemize}
		\item[(i)] There exists an edge $e=t_1t_2$ of $T$ such that for each $i \in [2]$, $G[A_{e,t_i}]$ is a $\xi$-nice graph with at least $\frac{m}{3}-\xi \geq \frac{m}{4}$ edges and with at most $\xi$ isolated vertices.
		\item[(ii)] There exists $t^* \in V(T)$ such that $\lvert E(G[X_{t^*}]) \rvert \geq \frac{m}{9}$.
		\item[(iii)] There exist $t^* \in V(T)$ and a partition $\{U_1,U_2\}$ of the set of components of $T-t^*$ such that for each $i \in [2]$, there are at least $\frac{2m}{9}$ edges of $G$ incident with $\bigcup_{C \in U_i}X_{V(C)}$.
	\end{itemize}

We first count the number of graphs satisfying (i).
For each $i \in [2]$, let $m_i = \lvert E(G[A_{e,t_i}]) \rvert$, so $G[A_{e,t_i}]$ is a $\xi$-nice graph with $m_i$ edges and with at most $\xi$ isolated vertices, and hence there are at most $\frac{c^{m_i}}{m_i^{3\xi^2+4}} \cdot (\xi+1)$ such graphs by the induction hypothesis.
Since each $G[A_{e,t_i}]$ has at most $\xi$ isolated vertices, $\lvert V(G[A_{e,t_i}]) \rvert \leq 2m_i+\xi$.
Therefore, there are at most $\frac{c^{m_1}}{m_1^{3\xi^2+4}}(\xi+1) \cdot \frac{c^{m_2}}{m_2^{3\xi^2+4}}(\xi+1) \cdot ((2m_1+\xi)(2m_2+\xi))^\xi \leq \frac{4^{3\xi^2+4}c^{m_1}}{m^{3\xi^2+4}} \cdot \frac{4^{3\xi^2+4}c^{m_2}}{m^{3\xi^2+4}} \cdot (2m)^{2\xi+2} \leq \frac{c^{m_1+m_2}}{m^{3\xi^2+4}} \cdot \frac{4^{6\xi^2+\xi+9}}{m^{3\xi^2+2-2\xi}} \leq \frac{c^m}{3m^{3\xi^2+4}}$ unlabelled $m$-edges $\xi$-nice graphs with no isolated vertex satisfying (i).

Second, we count the number of graphs satisfying (ii).
Let $m_0 = \lvert E(G[X_{t^*}]) \rvert$.
Then $G-E(G[X_{t^*}])$ is an $(m-m_0)$-edge $\xi$-nice graph with at most $\lvert X_{t^*} \rvert \leq \xi$ isolated vertices.
If $m_0=m$, then there are at most $m_0^{{\lvert X_{t^*} \rvert \choose 2}} \leq m^{\xi^2} \leq \frac{c^m}{3m^{3\xi^2+4}}$ such graphs, since $m>a$.
So we may assume that $m_0<m$.
By the induction hypothesis, there exist at most $\frac{c^{m-m_0}}{(m-m_0)^{3\xi^2+4}} \cdot (\xi+1)$ such $G-E(G[X_{t^*}])$.
Since $m_0 \geq \frac{m}{9}$ and $m>a$, there are at most $m_0^{{\lvert X_{t^*} \rvert \choose 2}} \cdot \frac{c^{m-m_0}}{(m-m_0)^{3\xi^2+4}} \cdot (\xi+1) \leq \frac{c^m}{m^{3\xi^2+4}} \cdot (\frac{m}{m-m_0})^{3\xi^2+4} \cdot c^{-m_0} \cdot m^{\xi^2}(\xi+1) \leq \frac{c^m}{m^{3\xi^2+4}} \cdot c^{-m/9} \cdot m^{4\xi^2+4}(\xi+1) \leq \frac{c^m}{3m^{3\xi^2+4}}$ $m$-edge unlabelled $\xi$-nice graphs with no isolated vertex satisfying (ii).

Now we count the number of graphs satisfying (iii) but not satisfying (ii).
Let $G_1 = G[X_{t^*} \cup \bigcup_{C \in U_1}X_{V(C)}]$.
Let $G_2 = G[X_{t^*} \cup \bigcup_{C \in U_2}X_{V(C)}] - E(G_1)$.
Note that $G_1$ and $G_2$ are edge-disjoint subgraphs of $G$.
In addition, the adhesion of $(T,\X)$ is at most $\xi$ and the torso at $t^*$ is obtained from a graph on at most $\xi$ vertices by attaching leaves adjacent to $X_{t^*}$, so there are at most $\xi^2$ edges of $G$ not belonging to $E(G_1) \cup E(G_2)$.
For $i \in [2]$, let $m_i$ be the number of edges of $G_i$, so $m_i \geq \frac{2m}{9}-\xi^2- \lvert E(G[X_{t^*}]) \rvert \geq \frac{m}{9}-\xi^2 \geq \frac{m}{10}$ and $m_1+m_2 \leq m$.
So for each $i \in [2]$, there exist at most $\xi^2+\lvert X_{t^*} \rvert \leq \xi^2+\xi$ isolated vertices in $G_i$, so it has at most $2m_i+\xi^2+\xi$ vertices.
By the induction hypothesis, there are at most $(\xi^2+\xi+1) \frac{c^{m_1}}{m_1^{3\xi^2+4}}$ possible unlabelled $G_1$, and there are at most $\sum_{j=1}^\xi{m_2 \choose j} j! \cdot (\xi^2+\xi+1) \frac{c^{m_2}}{m_2^{3\xi^2+4}} \leq m^{\xi+1} \frac{c^{m_2}}{m_2^{3\xi^2+4}}$ possible $G_2$ with at most $\xi$ vertices labelled by the elements in $X_{t^*}$ and with other vertices unlabelled.
	Therefore, there exist at most $(\xi^2+\xi+1) \frac{c^{m_1}}{m_1^{3\xi^2+4}} \cdot m^{\xi+1} \frac{c^{m_2}}{m_2^{3\xi^2+4}} \cdot ((2m_1+\xi^2+\xi)(2m_2+\xi^2+\xi))^{\xi^2} \leq \frac{c^{m_1+m_2}}{m^{3\xi^2+4} \cdot m^{3\xi^2+4}} \cdot {10}^{3\xi^2+4} \cdot {10}^{3\xi^2+4} \cdot 3\xi^2 \cdot m^{\xi+1} \cdot (9m^2)^{\xi^2} \leq \frac{c^{m_1+m_2}}{3m^{3\xi^2+4}} \cdot m^{-(3\xi^2+4)} \cdot m^{\xi+4} \cdot m^{2\xi^2} \leq \frac{c^m}{3m^{3\xi^2+4}}$ unlabelled $m$-edge $\xi$-nice graphs with no isolated vertex satisfying (iii).

Therefore, there are at most $3 \cdot \frac{c^m}{3m^{3\xi^2+4}} = \frac{c^m}{m^{3\xi^2+4}}$ unlabelled $m$-edge $\xi$-nice graphs with no isolated vertex.
\end{pf}

\bigskip

The following is a restatement of Theorem \ref{counting_basic_intro}.

\begin{theorem}
Let $d,h$ be positive integers with $d \geq 4$.
Then there exists a positive integer $c$ such that for every graph $H$ with maximum degree at most $d$ on at most $h$ vertices, and for every positive integer $m$, there are at most $c^m$ unlabelled $m$-edge $H$-immersion free graphs with no isolated vertex whose every maximal 2-edge-connected subgraph is $d$-edge-connected. 

Furthermore, there exists a positive integer $b$ such that for every graph $H$ with maximum degree at most $d$ on at most $h$ vertices and for every positive integer $n$, there are at most $b^n$ unlabelled $n$-vertex simple $H$-immersion free graphs whose every maximal 2-edge-connected subgraph is $d$-edge-connected.
\end{theorem}

\begin{pf}
Let $\xi = \xi_{\ref{mw_equiv}}(d,h)$, where $\xi_{\ref{mw_equiv}}$ is the integer $\xi$ mentioned in Corollary \ref{mw_equiv}.
Define $c=c_{\ref{counting_nice}}(2\xi+1)$, where $c_{\ref{counting_nice}}$ is the integer $c$ mentioned in Lemma \ref{counting_nice}.
Define $b=c^{7h+7}$.

Let $H$ be a graph with maximum degree at most $d$ on at most $h$ vertices.
By Corollary \ref{mw_equiv}, every $d$-edge-connected graph with no $H$-immersion is $2\xi$-nice.
Hence it is clear that every connected $H$-immersion free graph whose every maximal 2-edge-connected subgraph is $d$-edge-connected is $2\xi$-nice. 
So every $H$-immersion free graph whose every maximal 2-edge-connected subgraph is $d$-edge-connected is $(2\xi+1)$-nice.
Therefore, by Lemma \ref{counting_nice}, there are at most $c^m$ unlabelled $m$-edge $H$-immersion free graphs with no isolated vertex such that every maximal 2-edge-connected subgraph is $d$-edge-connected. 

Since every simple graph with no $H$-immersion has no $K_h$-immersion, by \cite[Theorem 1.4]{glw}, every simple $n$-vertex graph with no $K_h$-immersion has at most $(7h+6)n$ edges and has at most $n$ isolated vertices.
So there are at most $(n+1)c^{(7h+6)n} \leq b^n$ unlabelled $n$-vertex simple $H$-immersion free graphs whose every 2-edge-connected subgraph is $d$-edge-connected. 
\end{pf}

\end{document}